\documentclass[10pt]{article}
\usepackage{amsmath}
\usepackage{latexsym}
\usepackage{amssymb}
\usepackage{amsfonts}
\usepackage{amsthm}
\title{MACAULAY INVERSE SYSTEMS AND \\CARTAN-KAHLER THEOREM}
\author{J.-F. Pommaret \\ CERMICS, Ecole des Ponts ParisTech,\\ 6/8 Av. Blaise Pascal, 77455 Marne-la-Vall\'ee Cedex 02, France \\
E-mail: jean-francois.pommaret@wanadoo.fr, pommaret@cermics.enpc.fr \\
URL: http://cermics.enpc.fr/$\sim$pommaret/home.html }
\date{  }
\begin{document}
\maketitle

\noindent
{\bf ABSTRACT}\\

During the last months or so we had the opportunity to read two papers trying to relate the study of Macaulay (1916) {\it inverse systems} with the so-called Riquier (1910)-Janet (1920) {\it initial conditions} for the integration of linear analytic systems of partial differential equations. One paper has been written by F. Piras (1998) and the other by U. Oberst (2013), both papers being written in a rather algebraic style though using quite different techniques. It is however evident that the respective authors, though knowing the computational works of C. Riquier (1853-1929), M. Janet (1888-1983) and W. Gr\"{o}bner (1899-1980) done during the first half of the last century in a way not intrinsic at all, are not familiar with the formal theory of systems of ordinary or partial differential equations developped by D.C. Spencer (1912-2001) and coworkers around 1965 in an intrinsic way, in particular with its application to the study of {\it differential modules} in the framework of {\it algebraic analysis}. As a byproduct, the first purpose of this paper is to establish a close link between the work done by F. S. Macaulay (1862-1937) on inverse systems in 1916 and the well-known Cartan-K\"{a}hler theorem (1934). The second purpose is also to extend the work of Macaulay to the study of arbitrary linear systems with variable coefficients. The reader will notice how powerful and elegant is the use of the Spencer operator acting on sections in this general framework. However, we point out the fact that the literature on differential modules mostly only refers to a {\it complex analytic structure} on manifolds while the Spencer sequences have been created in order to study {\it any kind of structure} on manifolds defined by a Lie pseudogroup of transformations, {\it not just only complex analytic ones}. Many tricky explicit examples illustrate the paper, including the ones provided by the two authors quoted but in a quite different framework. \\

\noindent
{\bf 1 INTRODUCTION}\\

   With only a slight abuse of language, one can say that the birth of the {\it formal theory} of systems of ordinary differential (OD) or partial differential (PD) equations is coming from the combined work of C. Riquier ([22],1910) and M. Janet ([5],1920) along algebraic ideas brought by D. Hilbert in his study of sygygies. Roughly, one can say that the given OD or PD equations and all their derivatives may allow to compute a certain number of derivatives of the unknowns, called {\it principal}, from the other ones, called {\it parametric}, which can be chosen arbitrarily, on the condition of course that the resulting computation of this separation, also called "{\it cut}", should be unique or at least well defined. An (apparently) independent though similar approach has been followed later on by W. Gr\"{o}bner ([4], 1939). However, the successives approaches of Riquier, Janet and Gr\"{o}bner both suffer from the same lack of {\it intrinsicness} as they highly depend on the ordering of the $n$ independent variables and derivatives of the $m$ unknowns involved in a system of order $q$ ([13],[17]). In actual practice and as a summarizing comment, we may say that the central concepts in the successive works of Riquier, Macaulay ([9],$\S 38$, p 39), Janet and Gr\"{o}bner were indeed {\it formal integrability} and {\it involution}, almost fifty years before these concepts were properly defined and studied by Spencer and coworkers around 1965 ([13],[14],[25]).  \\
At the same time and as a way to generalize the situation to be found in the study of dynamical systems where a given {\it initial point} must be given in order 
to determine a unique trajectory passing through that point, a problem was raised by physicists working with many complicate field equations, namely that of the knowledge of the so-called "{\it degree of generality}" of solutions or at least a possibility to separate the many parametric derivatives into a certain number of "blocks" being described by arbitrary functions of a certain number of independent variables (It is now known from the published "letters on absolute parallelism" exchanged between E. Cartan and A. Einstein during the years 1929-1932 why Cartan only wrote one paper on the subject to Einstein in the language of PD equations, never quoting that it was a straight copy of the work done by Janet who suffered a lot about this situation as he told himself to the author of this paper while he was alife). As we shall see in the second section of this paper, the solution of this problem has first been described along the famous Cauchy-Kowaleski theorem and extended later on in 1934 along the Cartan-K\"{a}hler (CK) theorem ([3],[6]). The main problem, at least in our opinion, is that such a theorem is always presented in the framework of the exterior calculus of Cartan and thus totally separated from its formal origin which was essentially based on the involution assumption leading to the underlying Hilbert polynomial. In particular, the very specific type of systems met in the Spencer sequences cannot be imagined from the only use of the Janet sequences as we shall see.  \\

     Meanwhile, {\it commutative algebra}, namely the study of modules over rings, was facing a very subtle problem, the resolution of which led to the modern but difficult {\it homological algebra} with sequences and diagrams. Roughly, one can say that the problem was essentially to study properties of finitely generated modules not depending on the {\it presentation} of these modules by means of generators and relations. This very hard step is based on homological/cohomological methods like the so-called {\it extension} modules which cannot therefore be avoided ([1],[7],[15]) but are quite far from exterior calculus. Using now rings of {\it differential operators} instead of polynomial rings led to {\it differential modules} and to the challenge of adding the word {\it differential} in front of concepts of commutative algebra. Accordingly, not only one needs properties not depending on the presentation as we just explained but also properties not depending on the coordinate system as it becomes clear from any application to mathematical or engineering physics where tensors and exterior forms are always to be met like in the space-time formulation of electromagnetism. Unhappily, no one of the previous techniques for OD or PD equations could work.\\
         
     By chance, the intrinsic study of systems of OD or PD equations has been pioneered in a totally independent way by D. C. Spencer and collaborators after 1965 ([25]), as we already said, in order to relate differential properties of the PD equations to algebraic properties of their symbols, a technique superseding the {\it leading term} approach of Riquier, Macaulay, Janet or Gr\"{o}bner . Accordingly, it was another challenge to unify the {\it purely differential} approach of Spencer with the {\it purely algebraic} approach of commutative algebra, having in mind the necessity to use the previous homological algebraic results in this new framework. This sophisticated mixture of differential geometry and homological algebra, now called {\it algebraic analysis}, has been achieved after 1970 by M. Kashiwara [7] for the variable coefficients case.\\
     
    In a rough way, we shall prove in the third section of this paper that, if a {\it differential module} $M$ is defined over a differential field $K$ by a linear involutive system of OD or PD equations of any order $q$ with $n$ independent variables and $m$ unknowns, one can always find an isomorphic differential module defined by a linear involutive system in {\it Spencer form}, that is a first order involutive system not containing any zero order equation. Starting afresh with this new system, the CK {\it data} are made by a certain number of formal power series in $0$ variable (constants), $1$ variable, ..., up to $n$ variables, {\it the total number of such formal power series being equal to the number of unknowns}. Moreover, they allow to fully describe the dual {\it inverse system} $R=hom_K(M,K)$ as a left differential module for the Spencer operator acting on sections in the differential geometric framework. In the case of systems with constant or even variable coefficients as well, this result allows to exhibit a finite basis of $R$. Many explicit examples will illustrate this paper and must become test examples for using computer algebra without refering to Gr\"{o}bner bases. The reader will however notice that the concepts presented and the language of sections have rarely been used in mathematics (See [18],[19] for other details) and, up to our knowledge, have never been used in computer algebra or in mathematical physics (See [20] and [21] for very recent papers). \\

\noindent
{\bf 2  CARTAN-KHALER THEOREM REVISITED} \\

If $E$ is a vector bundle over the base manifold $X$ with projection $\pi$ and local coordinates $(x,y)=(x^i,y^k)$ projecting onto $x=(x^i)$ for $i=1,...,n$ and $k=1,...,m$, identifying a map with its graph, a (local) section $f:U\subset X \rightarrow E$ is such that $\pi\circ f =id$ on $U$ and we write $y^k=f^k(x)$ or simply $y=f(x)$. For any change of local coordinates $(x,y)\rightarrow (\bar{x}=\varphi(x),\bar{y}=A(x)y)$ on $E$, the change of section is $y=f(x)\rightarrow \bar{y}=\bar{f}(\bar{x})$ such that ${\bar{f}}^l(\varphi(x)\equiv A^l_k(x)f^k(x)$. The new vector bundle $E^*$ obtained by changing the {\it transition matrix} $A$ to its inverse $A^{-1}$ is called the {\it dual vector bundle} of $E$. Differentiating with respect to $x^i$ and using new coordinates $y^k_i$ in place of ${\partial}_if^k(x)$, we obtain ${\bar{y}}^l_r{\partial}_i{\varphi}^r(x)=A^l_k(x)y^k_i+{\partial}_iA^l_k(x)y^k$. Introducing a multi-index $\mu=({\mu}_1,...,{\mu}_n)$ with length $\mid \mu \mid={\mu}_1+...+{\mu}_n$ and prolonging the procedure up to order $q$, we may construct in this way, by patching coordinates, a vector bundle $J_q(E)$ over $X$, called the {\it jet bundle of order} $q$ with local coordinates $(x,y_q)=(x^i,y^k_{\mu})$ with $0\leq \mid\mu\mid \leq q$ and $y^k_0=y^k$. For a later use, we shall set $\mu+1_i=({\mu}_1,...,{\mu}_{i-1},{\mu}_i+1,{\mu}_{i+1},...,{\mu}_n)$ and define the operator $j_q:E \rightarrow J_q(E):f \rightarrow j_q(f)$ on sections by the local formula $j_q(f):(x)\rightarrow({\partial}_{\mu}f^k(x)\mid 0\leq \mid\mu\mid \leq q,k=1,...,m)$. Finally, a jet coordinate $y^k_{\mu}$ is said to be of {\it class} $i$ if ${\mu}_1=...={\mu}_{i-1}=0, {\mu}_i\neq 0$. \\

\noindent
{\bf DEFINITION 2.1}:  A {\it system} of PD equations of order $q$ on $E$ is a vector subbundle $R_q\subset J_q(E)$ locally defined by a constant rank system of linear equations for the jets of order $q$ of the form $ a^{\tau\mu}_k(x)y^k_{\mu}=0$. Its {\it first prolongation} $R_{q+1}\subset J_{q+1}(E)$ will be defined by the equations $ a^{\tau\mu}_k(x)y^k_{\mu}=0, a^{\tau\mu}_k(x)y^k_{\mu+1_i}+{\partial}_ia^{\tau\mu}_k(x)y^k_{\mu}=0$ which may not provide a system of constant rank as can easily be seen for $xy_x-y=0 \Rightarrow xy_{xx}=0$ where the rank drops at $x=0$.\\

The next definition will be crucial for our purpose.\\

\noindent
{\bf DEFINITION 2.2}: A system $R_q$ is said to be {\it formally integrable} if the $R_{q+r}$ are vector bundles $\forall r\geq 0$ (regularity condition) and no new equation of order $q+r$ can be obtained by prolonging the given PD equations more than $r$ times, $\forall r\geq 0$.\\

Finding an inrinsic test has been achieved by D.C. Spencer in 1965 along coordinate dependent lines sketched by Janet as early as in 1920 ([5]) and Gr\"{o}bner in 1940 ([5]), as we already said. The key ingredient, missing explicitly before the moderrn approach, is provided by the following definition.\\

\noindent
{\bf DEFINITION 2.3}: The family $g_{q+r}$ of vector spaces over $X$ defined by the purely linear equations $ a^{\tau\mu}_k(x)v^k_{\mu+\nu}=0$ for $ \mid\mu\mid= q, \mid\nu\mid =r $ is called the {\it symbol} at order $q+r$ and only depends on $g_q$.\\

The following procedure, {\it where one may have to change linearly the independent variables if necessary}, is the heart towards the next definition which is intrinsic even though it must be checked in a particular coordinate system called $\delta$-{\it regular} (See [13] and [14] for more details):\\

\noindent
$\bullet$ {\it Equations of class} $n$: Solve the maximum number ${\beta}^n_q$ of equations with respect to the jets of order $q$ and class $n$. Then call $(x^1,...,x^n)$ {\it multiplicative variables}.\\
\[  - - - - - - - - - - - - - - - -  \]
$\bullet$ {\it Equations of class} $i$: Solve the maximum number of {\it remaining} equations with respect to the jets of order $q$ and class $i$. Then call $(x^1,...,x^i)$ {\it multiplicative variables} and $(x^{i+1},...,x^n)$ {\it non-multiplicative variables}.\\
\[ - - - - - - - - - - - - - - - - - \]
$\bullet$ {\it Remaining equations equations of order} $\leq q-1$: Call $(x^1,...,x^n)$ {\it non-multiplicative variables}.\\

\noindent
{\bf DEFINITION 2.4}: A system of PD equations is said to be {\it involutive} if its first prolongation can be achieved by prolonging its equations only with respect to the corresponding multiplicative variables. The numbers ${\alpha}^i_q=m(q+n-i-1)!/((q-1)!(n-i)!)-{\beta}^i_q$ will be called {\it characters} and ${\alpha}^1_q\geq ... \geq {\alpha}^n_q $. For an involutive system, $(y^{{\beta}^n_q +1},...,y^m)$ can be given arbitrarily.  \\

Though the preceding description was known to Janet (he called it : "modules de formes en involution"), surprisingly he never used it explicitly. In any case, such a definition is far from being intrinsic and the hard step will be achieved from the Spencer cohomology that will also play an important part in the so-called {\it reduction to first order}, a result no so well known today as we shall see.\\

Let us consider $J_{q+1}(E)$ with jet coordinates $\{ y^l_{\lambda}\mid 0\leq \mid\lambda\mid\leq q+1\}$ and $J_1(J_q(E))$ with jet coordinates $\{z^k_{\mu},z^k_{\mu,i}\mid 0\leq \mid\mu\mid\leq q,i=1,...,n\}$. The canonical inclusion $J_{q+1}(E)\subset J_1(J_q(E))$ is described by the {\it two kinds} of equations:\\
\[   z^k_{\mu,i}-z^k_{\mu+1_i}=0 ,     \hspace{3cm}     0\leq \mid\mu\mid\leq q-1  \]
\[  z^k_{\mu+1_j,i}-z^k_{\mu+1_i,j}=0  , \hspace{3cm}   \mid\mu\mid=q-1  \]
or using the parametrization $z^k_{\mu,i}=y^k_{\mu+1_i}$ for $\mid\mu\mid=q$ with $z^k_{\mu}=y^k_{\mu}, \forall 0\leq \mid \mu \mid \leq q$. \\

Let $T$ be the tangent vector bundle of vector fields on $X$, $T^*$ be the cotangent vector bundle of 1-forms on $X$ and ${\wedge}^sT^*$ be the vector bundle of s-forms on $X$ with usual bases $\{dx^I=dx^{i_1}\wedge ... \wedge dx^{i_s}\}$ where we have set $I=(i_1< ... <i_s)$. Also, let $S_qT^*$ be the vector bundle of symmetric q-covariant tensors. Moreover, if  $\xi,\eta\in T$ are two vector fields on $X$, we may define their {\it bracket} $[\xi,\eta]\in T$ by the local formula $([\xi,\eta])^i(x)={\xi}^r(x){\partial}_r{\eta}^i(x)-{\eta}^s(x){\partial}_s{\xi}^i(x)$ leading to the {\it Jacobi identity} $[\xi,[\eta,\zeta]]+[\eta,[\zeta,\xi]]+[\zeta,[\xi,\eta]]=0, \forall \xi,\eta,\zeta \in T$. We have also the useful formula $[T(f)(\xi),T(f)(\eta)]=T(f)([\xi,\eta])$ where $T(f):T(X)\rightarrow T(Y)$ is the tangent mapping of a map $f:X\rightarrow Y$. Finally, we may introduce the {\it exterior derivative} $d:{\wedge}^rT^*\rightarrow {\wedge}^{r+1}T^*:\omega={\omega}_Idx^I \rightarrow d\omega={\partial}_i{\omega}_Idx^i\wedge dx^I$ with $d^2=d\circ d\equiv 0$ in the {\it Poincar\'{e} sequence}:\\
\[  {\wedge}^0T^* \stackrel{d}{\longrightarrow} {\wedge}^1T^* \stackrel{d}{\longrightarrow} {\wedge}^2T^* \stackrel{d}{\longrightarrow} ... \stackrel{d}{\longrightarrow} {\wedge}^nT^* \longrightarrow 0  \]

In a purely algebraic setting, one has ([13],[15]):  \\

\noindent
{\bf PROPOSITION 2.5}: There exists a map $\delta:{\wedge}^sT^*\otimes S_{q+1}T^*\otimes E\rightarrow {\wedge}^{s+1}T^*\otimes S_qT^*\otimes E$ which restricts to $\delta:{\wedge}^sT^*\otimes g_{q+1}\rightarrow {\wedge}^{s+1}T^*\otimes g_q$ and ${\delta}^2=\delta\circ\delta=0$.\\

{\it Proof}: Let us introduce the family of s-forms $\omega=\{ {\omega}^k_{\mu}=v^k_{\mu,i}dx^I\}$ and set $(\delta\omega)^k_{\mu}=dx^i\wedge{\omega}^k_{\mu+1_i}$. We obtain at once $({\delta}^2\omega)^k_{\mu}=dx^i\wedge dx^j\wedge{\omega}^k_{\mu+1_i+1_j}=0$.\\
\hspace*{12cm} Q.E.D.  \\

The kernel of each $\delta$ in the first case is equal to the image of the preceding $\delta$ but this may no longer be true in the restricted case and we set:\\

\noindent
{\bf DEFINITION 2.6}: We denote by $B^s_{q+r}(g_q)\subseteq Z^s_{q+r}(g_q)$ and $H^s_{q+r}(g_q)=Z^s_{q+r}(g_q)/B^s_{q+r}(g_q)$ respectively the coboundary space, cocycle space and cohomology space at ${\wedge}^sT^*\otimes g_{q+r}$ of the restricted $\delta$-sequence which only depend on $g_q$ and may not be vector bundles. The symbol $g_q$ is said to be s-{\it acyclic} if $H^1_{q+r}=...=H^s_{q+r}=0, \forall r\geq 0$, {\it involutive} if it is n-acyclic and {\it finite type} if $g_{q+r}=0$ becomes trivially involutive for r large enough. Finally, $S_qT^*\otimes E$ is involutive $\forall q\geq 0$ if we set $S_0T^*\otimes E=E$. \\

The preceding results will be used in proving the following technical result that will prove to be quite useful for our purpose ([13],[17]).\\

\noindent
{\bf PROPOSITION 2.7}: One has the isomorphisms $J_1(J_q(E))/J_{q+1}(E)\simeq T^*\otimes J_q(E)/\delta (S_{q+1}T^*\otimes E)\simeq T^*\otimes J_{q-1}(E)\oplus \delta(T^*\otimes S_qT^*\otimes E)$. More generally one may define the {\it Spencer bundles} $C_r(E)=C_{q,r}(E)={\wedge}^rT^*\otimes J_q(E)/\delta ({\wedge}^{r-1}T^*\otimes S_{q+1}T^*\otimes E)$ with an isomorphism $C_r(E)\simeq \delta ({\wedge}^rT^*\otimes S_qT^*\otimes E)\oplus {\wedge}^rT^*\otimes J_{q-1}(E)$. In particular one has $C_0(E)=J_q(E)$ while $C_n(E)={\wedge}^rT^*\otimes J_{q-1}(E)$.\\

{\it Proof}: The first commutative ad exact diagram:\\
\[  \begin{array}{rcccccl}
      & 0 &  & 0 &  & 0 &    \\
      &\downarrow &  &\downarrow &  &  \downarrow &   \\
0\rightarrow &S_{q+1}T^*\otimes E &\rightarrow &T^*\otimes J_q(E) &\rightarrow &C_1(E)&\rightarrow 0  \\
      & \downarrow &  &\downarrow &  &\parallel &   \\
0\rightarrow &J_{q+1}(E)&\rightarrow &J_1(J_q(E)) & \rightarrow &C_1(E) & \rightarrow 0  \\
     & \downarrow &   &\downarrow  &  & \downarrow &  \\
0\rightarrow & J_q(E) & = & J_q(E) &\rightarrow & 0  \\
  & \downarrow &   & \downarrow &  &   &  \\
  &   0   &   &   0  &   &   &
  \end{array}   \]
 \noindent 
shows that $C_1(E)\simeq T^*\otimes J_q(E)/S_{q+1}T^*\otimes E $. The general case finally depends on the following second commutative and exact diagram by using a (non-canonical) splitting of the right column:\\  
\[   \begin{array}{cccccl}
       0   &   &   0   &   &   0   &   \\
    \downarrow &  & \downarrow &  & \downarrow  &   \\
 {\wedge}^{r-1}T^*\otimes S_{q+1}T^*\otimes E  & \stackrel{\delta}{ \rightarrow} & {\wedge}^rT^*\otimes S_qT^*\otimes E &
 \stackrel{\delta}{\rightarrow} & \delta ( {\wedge}^rT^*\otimes S_qT^*\otimes E) &  \rightarrow 0  \\
      \parallel  &  &  \downarrow  &  & \downarrow  &  \\
 {\wedge}^{r-1}T^*\otimes S_{q+1}T^*\otimes E &\rightarrow  &{\wedge}^rT^*\otimes J_q(E) & \rightarrow &C_r(E) &
 \rightarrow 0  \\
      \downarrow  &  & \downarrow &  & \downarrow &   \\
        0  &  \rightarrow &  {\wedge}^rT^*\otimes J_{q-1}(E) &  = & {\wedge}^r T^*\otimes J_{q-1}(E) & \rightarrow 0  \\
         &   &  \downarrow &  &  \downarrow  &  \\
         &   &  0  &  &  0  &  
       \end{array}     \]
When $r=n$, the equality $\delta ({\wedge}^{n-1}T^*\otimes S_{q+1}t^*)= {\wedge}^nT^*\otimes  S_qT^*$ gives the last result. \\
\hspace*{12cm} Q.E.D.  \\
    
These absolutely non-trivial results can be restricted to the systems and symbols. Accordingly, the inclusion $R_{q+1}\subset J_1(R_q)$ can be considered as a new first order system over $R_q$, called {\it first order reduction} or {\it Spencer form}. One obtains ([13],[14],[15]):\\
 
 \noindent
{\bf PROPOSITION 2.8}: The first order reduction is formally integrable (involutive) whenever $R_q$ is formally integrable (involutive). In that case, the reduction has no longer any zero order equation.\\
  
 Having in mind control theory, we have therefore set up the problem of "state", even for systems which are not of finite type and {\it it just remains to modify the Spencer form} in order to generalize the Kalman form to PD equations. Here is the procedure that must be followed in the case of a first order involutive system with no zero order equation, for example like the one we just obtained.\\

\noindent 
$\bullet$  Look at the equations of class n solved with respect to $y^1_n,...,y^{\beta}_n$.\\
$\bullet$  Use integrations by part like:\\
\[ y^1_n-a(x)y^{\beta +1}_n=d_n(y^1-a(x)y^{\beta +1})+{\partial}_na(x)y^{\beta +1}={\bar{y}}^1_n+{\partial}_na(x)y^{\beta +1}  \]
$\bullet$  Modify $y^1,...,y^{\beta} $ to ${\bar{y}}^1,...,{\bar{y}}^{\beta}$ in order to "{\it absorb}" the various $y^{\beta +1}_n,...,y^m_n$ {\it only appearing in the equations of class} n.  \\

\noindent
{\bf PROPOSITION 2.9}: The new equations of class n only contain $y^{\beta +1}_i,...,y^m_i$ with $0\leq i\leq n-1$ while the equations of class $1,...,n-1$ no more contain $y^{\beta+1},...,y^m$ and their jets.\\

{\it Proof}: The first assertion comes from the absorption procedure. Now, if $y^m$ or $y^m_i$ should appear in an equation of class $\leq n-1$, prolonging this equation with respect to the non-multiplicative variable $x^n$ should bring $y^m_n$ or $y^m_{in}$ and (here involution is essential) we should get a linear combination of equations of various classes prolonged with respect to $x^1,...,x^{n-1}$ {\it only}, but this is impossible.\\
\hspace*{12cm} Q.E.D.  \\

A similar proof provides at once (See later on for the definition):\\

\noindent
{\bf COROLLARY 2.10}: Any torsion element, if it exists, only depends on ${\bar{y}}^1,...,{\bar{y}}^{\beta}$.\\

We are now in position to revisit Gr\"{o}bner bases with critical eyes (See [17] for more details).\\

\noindent
{\bf EXAMPLE 2.11}: Let $P_1=({\chi}_3)^2,P_2={\chi}_2{\chi}_3-({\chi}_1)^2,P_3=({\chi}_2)^2$ be three polynomials generating the ideal $\mathfrak{a}=(P_1,P_2,P_3)\subset \mathbb{Q}[{\chi}_1,{\chi}_2,{\chi}_3]$. The corresponding system $R_2$ defined by the three PD equations\\
\[     y_{33}=0, \hspace{5mm}y_{23}-y_{11}=0,\hspace{5mm} y_{22}=0  \]
is homogeneous and thus automatically formally integrable but $g_2$ is not involutive though finite type because $g_4=0$ (Exercise). Elementary computations of ranks of matrices shows that the $\delta$-map:\\
\[    0\rightarrow  {\wedge}^2T^*\otimes g_3  \stackrel{\delta}{\longrightarrow} {\wedge}^3T^*\otimes g_2 \rightarrow 0  \]
is an isomorphism and thus $g_3$ is 2-acyclic, a crucial intrinsic property [13,15,25] totally absent from any "old" work. Now, denoting the {\it initial} of a 
polynomial by $in(   )$ while choosing the ordering ${\chi}_3>{\chi}_1>{\chi}_2$, we obtain:  \\
\[  in(P_1)=({\chi}_3)^2, in(-P_2)=({\chi}_1)^2, in(P_3)=({\chi}_2)^2  \]
and $\{P_1,P_2,P_3\}$ is a Gr\" {o}bner basis. However, choosing the ordering ${\chi}_3>{\chi}_2>{\chi}_1$, we have now:\\
\[  in(P_1)=({\chi}_3)^2, in(P_2)={\chi}_2{\chi}_3, in(P_3)=({\chi}_2)^2  \]
and $\{P_1,P_2,P_3\}$ is not a Gr\"{o}bner basis because $ y_{112}=0, y_{113}=0$ AND $y_{1111}=0 $. Accordingly, a Gr\" {o}bner basis could be $\{P_1,P_2,P_3,P_4=({\chi}_1)^2{\chi}_2,P_5=({\chi}_1)^2{\chi}_3,P_6=({\chi}_1)^4\}$ (!!!).\\

The main use of involution is to construct differential sequences made up by successive {\it compatibility conditions} (CC). In particular, when $R_q$ is involutive, the linear differential operator ${\cal{D}}:E\stackrel{j_q}{\rightarrow} J_q(E)\stackrel{\Phi}{\rightarrow} J_q(E)/R_q=F_0$ of order $q$ with space of solutions $\Theta\subset E$ is said to be {\it involutive} and one has the canonical {\it linear Janet sequence} ([19], p 144):\\
\[  0 \longrightarrow  \Theta \longrightarrow T \stackrel{\cal{D}}{\longrightarrow} F_0 \stackrel{{\cal{D}}_1}{\longrightarrow}F_1 \stackrel{{\cal{D}}_2}{\longrightarrow} ... \stackrel{{\cal{D}}_n}{\longrightarrow} F_n \longrightarrow 0   \]
where each other operator is first order involutive and generates the CC of the preceding one with $F_r={\wedge}^rT^*\otimes J_q(E)/({\wedge}^rT^*\otimes R_q+\delta ({\wedge}^{r-1}T^*\otimes S_{q+1}T^*\otimes E))$. As the Janet sequence can be "cut at any place", that is can also be constructed anew from any intermediate operator, {\it the numbering of the Janet bundles has nothing to do with that of the Poincar\'{e} sequence for the exterior derivative}, contrary to what many physicists  still believe. Moreover, the fiber dimension of the Janet bundles can be computed at once inductively from the board of multiplicative and non-multiplicative variables that can be exhibited for $\cal{D}$ by working out the board for ${\cal{D}}_1$ and so on. For this, the number of rows of this new board is the number of dots appearing in the initial board while the number $nb(i)$ of dots in the column $i$ just indicates the number of CC of class $i$ for $i=1, ... ,n$ with $nb(i) < nb(j), \forall i<j$ and we have therefore:  \\

\noindent
{\bf THEOREM 2.12}: The successive first order operators ${\cal{D}}_1, ... , {\cal{D}}_n$ are {\it automatically} in reduced Spencer form. \\

\noindent
{\bf DEFINITION 2.13}: The Janet sequence is said to be {\it locally exact at} $F_r$ if any local section of $F_r$ killed by ${\cal{D}}_{r+1}$ is the image by ${\cal{D}}_r$ of a local section of $F_{r-1}$. It is called {\it locally exact} if it is locally exact at each $F_r$ for $0\leq r \leq n$. The Poincar\'{e} sequence is locally exact, that is a closed form is locally an exact form but counterexamples may exist ([20], p 373).\\

\noindent
{\bf EXAMPLE 2.14}: ([9],$\S 38$, p 40 where one can find the first intuition of formal integrability) The primary ideal $\mathfrak{q}=(({\chi}_1)^2, {\chi}_1{\chi}_3-{\chi}_2)$ provides the system $y_{11}=0, y_{13}-y_2=0$ which is neither formally integrable nor involutive. Indeed, we get $d_3y_{11}-d_1(y_{13}-y_2)=y_{12}$ and $d_3y_{12}-d_2(y_{13}-y_2)=y_{22}$, that is to say {\it each first and second} prolongation does bring a new second order PD equation. Considering the new system $y_{22}=0, y_{12}=0, y_{13}-y_2=0, y_{11}=0$, the question is to decide whether this system is involutive or not. One could use Janet/Gr\"{o}bner algorithm but with no insight towards involution. In such a simple situation, as there is no PD equation of class $3$, the evident permutation of coordinates $(1,2,3)\rightarrow (3,2,1)$ provides the following involutive second order system with one equation of class $3$, $2$ equations of class $2$ and $1$ equation of clas $1$:   \\
\[  \left\{  \begin{array}{lcl}
{\Phi}^4 \equiv y_{33}  & = & 0  \\
{\Phi}^3 \equiv y_{23} & = & 0  \\
{\Phi}^2 \equiv y_{22} & = &  0  \\
{\Phi}^1 \equiv y_{13}-y_2 & = &  0 
\end{array}
\right. \fbox{$\begin{array}{lll}
1 & 2 & 3 \\
1 & 2 & \bullet \\
1 & 2 & \bullet \\
1 & \bullet & \bullet
\end{array}$}  \]
We have ${\alpha}^3_2=0,{\alpha}^2_2=0,{\alpha}^1_2=2$ and the corresponding CC system is easily seen to be the following involutive first order system in reduced Spencer form:  \\
\[  \left\{  \begin{array}{lcl}
{\Psi}^4 \equiv  d_3{\Phi}^3-d_2{\Phi}^4  & = & 0  \\
{\Psi}^3 \equiv  d_3{\Phi}^2-d_2{\Phi}^3  & = & 0  \\
{\Psi}^2 \equiv  d_3{\Phi}^1-d_1{\Phi}^4+{\Phi}^3 & = &  0  \\
{\Psi}^1 \equiv  d_2{\Phi}^1-d_1{\Phi} ^3+ {\Phi}^2   & = &  0 
\end{array}
\right. \fbox{$\begin{array}{lll}
1 & 2 & 3  \\
1 & 2 & 3  \\
1 & 2 & 3  \\
1 & 2 & \bullet
\end{array}$}  \]
The final CC system is the involutive first order system in reduced Spencer form:  \\
\[  \left\{  \begin{array}{lcl}
\Omega \equiv  d_3{\Psi}^1-d_2{\Psi}^2+d_1{\Psi}^4-{\Psi}^3   & = & 0  
\end{array}
\right. \fbox{$\begin{array}{lll}
1 & 2 & 3 
\end{array}$}  \]
We get therefore the Janet sequence:    
\[    0 \longrightarrow  \Theta \longrightarrow 1 \longrightarrow 4 \longrightarrow 4 \longrightarrow 1  \longrightarrow  0    \]
and check that the Euler-Poincar\'{e} characteristic, that is the alternate sum of dimensions of the Janet bundles, is $1-4+4-1=0$.  \\ 

Equivalently, we have the involutive {\it first Spencer operator} $D_1:C_0=R_q\stackrel{j_1}{\rightarrow}J_1(R_q)\rightarrow J_1(R_q)/R_{q+1}\simeq T^*\otimes R_q/\delta (g_{q+1})=C_1$ of order one induced by the {\it Spencer operator} $D:R_{q+1}\rightarrow T^*\otimes R_q:{\xi}_{q+1} \rightarrow j_1({\xi}_q)-{\xi}_{q+1}$ which is well defined because both $J_{q+1}(E)$ and $T^*\otimes J_q(E)$ may be considered as sub-bundles of $J_1(J_q(E))$. Introducing the {\it Spencer bundles} $C_r={\wedge}^rT^*\otimes R_q/{\delta}({\wedge}^{r-1}T^*\otimes g_{q+1})$, the first order involutive ($r+1$)-{\it Spencer operator} $D_{r+1}:C_r\rightarrow C_{r+1}$ is induced by $D:{\wedge}^rT^*\otimes R_{q+1}\rightarrow {\wedge}^{r+1}T^*\otimes R_q:\alpha\otimes {\xi}_{q+1}\rightarrow d\alpha\otimes {\xi}_q+(-1)^r\alpha\wedge D{\xi}_{q+1}$ and we obtain the canonical {\it linear Spencer sequence} ([14], p 150):\\
\[    0 \longrightarrow \Theta \stackrel{j_q}{\longrightarrow} C_0 \stackrel{D_1}{\longrightarrow} C_1 \stackrel{D_2}{\longrightarrow} C_2 \stackrel{D_3}{\longrightarrow} ... \stackrel{D_n}{\longrightarrow} C_n\longrightarrow 0  \]
as the canonical Janet sequence for the first order involutive system $R_{q+1}\subset J_1(R_q)$.\\

The canonical Janet sequence and the canonical Spencer sequence can be connected by a commutative diagram where the Spencer sequence is induced by the locally exact central horizontal sequence which is at the same time the Janet sequence for $j_q$ and the Spencer sequence for $J_{q+1}(E)\subset J_1(J_q(E))$ ([14], p 153) but this result will not be used in this paper (See [16],[19],[20],[21] for recent papers providing more details on applications of these results to engineering and mathematical physics, in particular continuum mechanics, gauge theory and general relativity).  \\

   For an involutive system of order $q$ in solved form, we shall use to denote by $y_{pri}$ the {\it principal jet coordinates}, namely the leading terms of the solved equations in the sense of involution. Accordingly, any formal derivative of a principal jet coordinate is again a principal jet coordinate and the remaining jet coordinates are the parametric jet coordinates denoted by $y_{par}$. We shall use a "trick" in order to study the remaining jet coordinates called {\it parametric jet coordinates} and denoted by $y_{par}$. Indeed, the symbol of $j_q$ is the zero symbol and is thus trivially involutive at any order $q$. Accordingly, if we introduce the {\it multiplicative variables} $x^1,...,x^i$ for the parametric jets of order $q$ and class $i$, the formal derivative or a parametric jet of strict order $q$ and class $i$ by one of its multiplicative variables is uniquely obtained and cannot be a principal jet of order $q+1$ which is coming from a uniquely defined principal jet of order $q$ and class $i$. We have thus obtained the following technical Proposition which is very useful in actual practice: \\
   
\noindent
{\bf PROPOSITION 2.15}: The principal and parametric jets of strict order $q$ of an involutive system of order $q$ have the same Janet board if we extend it to all the classes that may exist for both sets, in particular the respective empty classes.   \\
   
   Paying attention to the specific situation of the symbol of order $q$, the following technical lemmas are straightforward consequences of the definition of an involutive system and allow to construct all the possible sets of principal or parametric jet coordinates when $m,n$ and $q$ are given (See [13] p 123-125 for more details). \\
   
 \noindent
 {\bf LEMMA 2.16}: If $y^k_{\mu}\in y_{pri}$ and $y^l_{\nu}\in y_{par}$ appear in the same equation of class $i$ in solved form, then $\nu$ is of class $\leq i$ and $l>k$ when $\nu$ is also of class $i$.  \\
 
 \noindent
 {\bf LEMMA 2.17}: If $y^k_{\mu}$ is a principal jet coordinate of strict order $q$, that is $\mid \mu \mid = q$ with ${\mu}_1=0, ..., {\mu}_{i-1}=0, {\mu}_i> 0$, then $\forall j>i$, $y^k_{\mu-1_i+1_j}$ is a principal jet coordinate and this notation has a meaning because ${\mu}_i>0$. \\
 
 \noindent
{\bf LEMMA 2.18}: If there exists an equation of class $i$, there exists also an equation of class $i+1$. Accordingly, the classes of the solved equations of an involutive symbol are an increasing chain of consecutive integers ending at $n$.  \\

\noindent
{\bf LEMMA 2.19}: The indices ${\mu}_i$ of the principal jet coordinates of strict order $q$ and class $i$ are an increasing chain of consecutive integers starting from $1$.  \\

Combining the preceding lemmas, we obtain:  \\

\noindent
{\bf PROPOSITION 2.20}: If $q=1$, one has $0\leq {\beta}^1_1\leq ... \leq {\beta}^n_1\leq m $ in a coherent way with the relations ${\alpha}^i_1=m-{\beta}^i_1,\forall i=1,...,n$ and $ (m-{\beta}^n_1)+({\beta}^n_1-{\beta}^{n-1}_1)+ ... +({\beta}^1_1-0)=m$. The system made by the equations of class $1$+ ... + class $i$ is involutive over $K[d_1,...,d_i]$ for $i=1,...,n$.\\

Using the Janet board and the definition of involutivity, we get $dim(g_{q+r})={\sum}_{i=1}^n\frac{(r+i-1)!}{r!(i-1)!}{\alpha}^i_q$ and thus $dim(R_{q+r})=dim(R_{q-1})+{\sum}_{i=1}^n\frac{(r+i)!}{r!i!}{\alpha}^i_q$. In the case of analytic systems, the following theorem providing the CK {\it data} is well known though its link with involution is rarely quoted because it is usually presented within the framework of exterior calculus ([3],[7]):  \\

\noindent
{\bf THEOREM 2.21} (Cartan-K\"{a}hler): If $R_q\subset J_q(E)$ is a linear involutive and analytic system of order $q$ on $E$, there exists one analytic solution $y^k=f^k(x)$ and only one such that:  \\
1) $(x_0,{\partial}_{\mu}f^k(x_0))$ with $0 \leq \mid \mu \mid\leq q-1$ is a point of $R_{q-1}={\pi}^q_{q-1}(R_q)\subset J_{q-1}(E)$.  \\
2) For $i=1,...,n$ the ${\alpha}^i_q$ parametric derivatives ${\partial}_{\mu}f^k(x)$ of class $i$ are equal for $x^{i+1}=x^{i+1}_0,...,x^n=x^n_0$ to ${\alpha}^i_q$ given analytic functions of $x^1,...,x^i$.  \\
 
 The monomorphism $0\rightarrow J_{q+1}(E) \rightarrow J_1(J_q(E))$ allows to identify $R_{q+1} $ with its image ${\hat{R}}_1$ in $ J_1(R_q)$ and we just need to set $R_q=\hat{E}$ in order to obtain the first order system (Spencer form) ${\hat{R}}_1\subset J_1(\hat{E})$ which is also involutive and analytic while ${\pi}^1_0:{\hat{R}}_1\rightarrow \hat{E}$ is an epimorphism. Studying the respective symbols, we may identify $g_{q+r}$ and ${\hat{g}}_r$ while ${\hat{g}}_1$ is involutive. Looking at the Janet board of multiplicative variables we have:
 \[   {\hat{\alpha}}^i_1={\alpha}^i_q+ ... + {\alpha}^n_q={\alpha}^i_{q+1} \Rightarrow {\alpha}^i_q={\hat{\alpha}}^i_1-{\hat{\alpha}}^{i+1}_1={\hat{\beta}}^{i+1}_1-{\hat{\beta}}^i_1  \]
and obtain:  \\

\noindent
{\bf COROLLARY 2.22}: If $R_1\subset J_1(E)$ is a first order linear involutive and analytic system such that ${\pi}^1_0:R_1 \rightarrow E$ is an epimorphism, then there exists one analytic solution $y^k=f^k(x)$ and only one, such that:  \\
1) $f^1(x),...,f^{{\beta}^1_1}(x)$ are equal to ${\beta}^1_1$ given constants when $x=x_0$.  \\
2) $f^{{\beta}^i_1+1}(x),..., f^{{\beta}^{i+1}_1}(x)$ are equal to ${\beta}^{i+1}_1-{\beta}^i_1$ given analytic functions of $x^1,...,x^i$ when $x^{i+1}=x^{i+1}_0, ...,x^n=x^n_0$.  \\
3) $f^{{\beta}^n_1+1}(x),...,f^m(x)$ are $m-{\beta}^n_1$ given analytic functions of $x^1,...,x^n$.  \\

\noindent
{\it Proof}: The analytic proof of the corollary uses inductively the Cauchy-Kowaleski theorem which is the particular case described by the conditions ${\beta}^1_1=0, ..., {\beta}^{n-1}_1=0,{\beta}^n_1=m$ leading to the existence of one analytic solution $y^k=f^k(x)$ and only one such that $f^1(x),...,f^m(x)$ are equal to $m$ given analytic functions of $x^1,...,x^{n-1}$ when $x^n=x^n_0$. As such a proof is quite technical, we refer the reader to ([13],p 159-163) for the details. For the reader not familiar with involution, we shall nevertheless explain in a formal way why the CK data of the corollary coincide with the CK data of the theorem, {\it a result not evident at first sight}, even on the next elementary illustrative examples.    \\
Indeed, according to the CK theorem, we get from $1)$ the parametric jets $y^1,...,y^m$ with no multiplicative variable and from $2)$ the ${\alpha}^n_1=m-{\beta}^n_1$ parametric jets $y^{{\beta}^n_1+1}_n,...,y^m_n$ of class $n$ with multiplicative variables $x^1,...,x^n$, ... , ${\alpha}^i_1=m-{\beta}^i_1$ parametic jets $y^{{\beta}^i_1+1}_i, ... ,y^m_i$ of class $i$ with multiplicative variables $x^1,...,x^i$, ... , and ${\alpha}^1_1$ parametric jets $y^{{\beta}^1_1},...,y^m_1$ of class $1$ with multiplicative variable $x^1$ only. Collecting the $y^k$ and their jets successively from the class $n$ down to the class $1$, we may start with $y^{{\beta}^n_1+1}_n,...,y^m_n$ with multiplicative variables $x^1,...,x^n$, ... , then $y^{{\beta}^n_1+1}_i,...,y^m_i$ with multiplicative variables $x^1,...,x^i$ because ${\beta}^n_1\geq {\beta}^i_1$, ... , then $y^{{\beta}^n_1+1}_1,...,y^m_1$ with multiplicative variable $x^1$ and finally $y^{{\beta}^n_1+1},...,y^m$ with no multiplicative variable or, {\it equivalently}, $y^{{\beta}^n_1+1},...,y^m$ with multiplicative variables $x^1,...,x^n$. We just need to repeat this procedure with $n-1$ in place of $n$ and ${\beta}^n_1$ in place of $m$, so on down till $1$ and ${\beta}^1_1$ in order to obtain the corollary. It is important to notice that only the use of the Spencer form can bring the same total number of such formal power series of $n$ down to $0$ independent variables as the number of unknowns, {\it even for nonlinear systems or linear systems with variable coefficients}.  \\
\hspace*{12cm}  Q.E.D.  \\ 
  
\noindent
{\bf EXAMPLE 2.23}: With $n=2, m=2, q=1, k=\mathbb{Q}$, let us consider the following first order involutive sytem in (reduced) Spencer form:  \\
\[  \left\{  \begin{array}{lcl}
{\Phi}^3\equiv y^2_2-y^2_1 & = & 0 \\
{\Phi}^2\equiv y^1_2 & = & 0 \\
{\Phi}^1\equiv y^1_1 & = &  0 
\end{array}
\right. \fbox{$\begin{array}{ll}
1 & 2  \\
1 & 2 \\
1 & \bullet
\end{array} $}  \]
with corresponding board of principal/parametric jets up to order $1$: \\
\[  \left\{  \begin{array}{l}
   y^1_2,y^2_2\\
y^1_1 \\
{  }
\end{array}
\right. \fbox{$\begin{array}{cc}
1 & 2 \\
1 & \bullet\\
\bullet & \bullet 
\end{array} $}  
\begin{array}{r}
  \\
y^2_1 \\
y^1,y^2
\end{array}      \]
because the principal class $2$ is full.   \\
The previous equations and their first prolongations are described by the board (Compare to [9]): \\
\[  \hspace*{14mm}  \begin{array}{r|cc|cccc|cccccc|c}
  & y^1 & y^2 &y^1_1 & y^1_2 & y^2_1 & y^2_2 & y^1_{11} & y^1_{12} & y^1_{22} & y^2_{11} & y^2_{12} & y^2_{22} & ...  \\
\hline
{\Phi}^3 & 0 & 0 & 0 & 0 & -1 & 1 & 0 & 0 & 0 & 0 & 0 & 0 & ...   \\
{\Phi}^2 & 0 & 0 & 0 & 1 & 0 & 0 & 0 & 0 & 0 & 0 & 0 & 0 & ...  \\
{\Phi}^1 & 0 & 0 & 1 & 0 & 0 & 0 & 0 & 0 & 0 & 0 & 0 & 0 & ...  \\
d_1{\Phi}^3 & 0 & 0 & 0 & 0 & 0 & 0 & 0 & 0 & 0& -1& 1 & 0 & .. 
\end{array}       \]
Accordingly, a basis of sections may start like in the board (Compare to [9]):  \\
\[  \begin{array}{lcl|cc|cccc|cccccc|c}
   &  &  &y^1 & y^2 & y^1_1 & y^1_2 & y^2_1 & y^2_2 &y^1_{11} & y^1_{12} & y^1_{22} &y^2_{11} & y^2_{12} & y^2_{22} & ...  \\
  \hline
f^1_0   &   \rightarrow    &  E^0_1       &  1 & 0 & 0 & 0 & 0 & 0& 0& 0 & 0 &0 & 0 & 0 & ... \\
f^2_0   &   \rightarrow    &E^0_2         &  0 & 1 & 0 & 0 & 0 & 0 & 0 & 0& 0& 0 & 0 & 0 & ... \\
f^2_1    &  \rightarrow    &  E^1_2       &  0 & 0 & 0 & 0 & 1 & 1 & 0 & 0 & 0 & 0 & 0 & 0 &  ... \\
f^2_{11} &  \rightarrow  &  E^{11}_2  &  0 & 0 & 0 & 0 & 0 & 0 & 0 & 0 & 0 & 1 & 1 & 1 & ...
\end{array}    \]
As the system is homogeneous, we have for example $E^{11}_2\equiv a^{11}_2+a^{12}_2+a^{22}_2=0$ and so on. We obtain therefore an infinite number of modular equations through this way to write down sections.\\
The CK data are $\{f^1(0,0), f^2(x^1,0)\}$ in agrement with the above result and we have:   \\
\[  \begin{array}{rcl}
 f^1(x^1,x^2) &= & f^1(0,0) \\
 f^2(x^1,x^2) & = & f^2(x^1,0)+{\partial}_2f^2(x^1,0)x^2+ ...\\
                        & = &  f^2(x^1,0)+{\partial}_1f^2(x^1,0)x^2+ ...  
 \end{array}   \]
In the present situation, we notice that $f^2(x^1,x^2)=f(x^1+x^2)=f^2(x^1+x^2,0)$ by using the first equation. \\
With ${\Phi}^3\equiv y^2_2-y^2_1=0, {\Phi}^2\equiv y^1_2-ay^2=0, {\Phi}^1\equiv y^1_1-ay^2=0$ where $a$ is a constant parameter, we let the reader check that this system is involutive for any value of $a$ and may not be homogeneous but the corresponding module $M$ is $1$-pure if and only if $a\neq 0$.  \\
With now $n=2, m=3, q=1, K=k=\mathbb{Q}$, we could finally add a third unknown $y^3$ and consider the following new first order involutive system in reduced Spencer form where $a$ is an arbitrary constant parameter:  \\
\[  \left\{  \begin{array}{lcl}
y^2_2-y^2_1- a y^3 & = & 0 \\
y^1_2 & = & 0 \\
y^1_1& = & 0
\end{array}
\right. \fbox{$\begin{array}{lll}
1 & 2  \\
1 & 2 \\
1 & \bullet 
\end{array} $}  \]
At order $1$ we have $pri=\{ y^1_1,y^1_2,y^2_2\}$ and $par=\{ y^1,y^2, y^3, y^2_1, y^3_1, y^3_2\}$ with corresponding board of principal/parametric jets up to order $1$:  \\
\[  \left\{  \begin{array}{l}
y^1_2,y^2_2 \\
y^1_1 \\
 {  }
\end{array}
\right. \fbox{$\begin{array}{ll}
1 & 2  \\
1 & \bullet  \\
\bullet & \bullet 
\end{array} $}  
\begin{array}{r}
y^3_2  \\
y^2_1,y^3_1  \\
y^1,y^2,y^3
\end{array}    \]
We obtain at once ${\beta}^2_1=2> {\beta}^1_1=1$ and the characters ${\alpha}^2_1=3-2=1<{\alpha}^1_1=3-1=2$, that is $m-{\beta}^2_1=3-2=1,{\beta}^2_1-{\beta}^1_1=2-1=1, {\beta}^1_1-0=1-0=1$. We have in particular $y^3_2$ with multiplicative variables $x^1,x^2$, then $y^3_1$ with multiplicative variable $x^1$ only and finally $y^3$ with no multiplicative variable or, {\it equivalently}, $y^3$ with multiplicative variables $x^1,x^2$. Similarly, we have $y^2_1$ with multiplicative variable $x^1$ only and $y^2$ with no multiplicative variable or, {\it equivalently}, $y^2$ with multiplicative variable $x^1$ to which we have finally to add $y^1$ with no multiplicative variable. It follows that the CK data are $\{ f^1(0,0), f^2(x^1,0), f^3(x^1,x^2)\}$ and we have indeed:  \\ 
\[ \begin{array}{rcl}
 f^2(x^1,x^2) & = & f^2(x^1,0)+{\partial}_2f^2(x^1,0)x^2+{\partial}_{22}f^2(x^1,0)\frac{(x^2)^2}{2}+ ...  \\
                        & = & f^2(x^1,0)+({\partial}_1f^2(x^1,0)+ a f^3(x^1,0))x^2  \\
                         &  &  +({\partial}_{11}f^2(x^1,0)+ a {\partial}_1f^3(x^1,0)+ a {\partial}_2f^3(x^1,0))\frac{(x^2)^2}{2}+ ... 
\end{array}    \]
Finally, we have ${\partial}_2f^1=0$ and ${\partial}_1f^1=0$. Accordingly, we have $f^1(x^1,x^2)=cst= f^1(0,0) $ as a way to obtain $f^1(x^1,x^2)$ in this particular situation. Though the CK data do not depend on $a$, the underlying differential module $M$ highly depends on $a$. Indeed, when $a\neq 0$ the torsion module $t(M)$ is generated by $z=y^1$ satisfying $z_2=0, z_1=0$ and we have the purity filtration:  \\
\[   0=t_2(M)\subset t_1(M) = t_0(M) = t(M)\subset M   \]
while, when $a=0$, the torsion module $t(M)$ is generated by $z'=y^1$ and $z"=y^2$ with both $z'_2=0$ and $z'_1=0$ but $z"_2-z"_1=0$ only, a result leading to the different purity filtration:  \\
\[      0=t_2(M)\subset t_1(M) \subset t_0(M)=t(M) \subset M  \]

\noindent
{\bf EXAMPLE 2.24}: With $n=2, m=1, q=8, k=\mathbb{Q}$, let us revisit the example of ([12],p 93). Using the multi-index notation, let us consider the system of seventh order $y_{(3,4)}=0, y_{(5,2)}=0$. This homogeneous system is of course formally integrable but is far from being involutive. Using the change of variables $x^1\rightarrow x^1+x^2, x^2 \rightarrow x^2$we get $\delta$-regular coordinates but the system/symbol is not involutive and we need one prolongation in order to obtain the following eighth order involutive system:  \\
\[  \left\{  \begin{array}{lcl}
y_{(0,8)}-15y_{(4,4)}-24y_{(5,3)}-10y_{6,2)} & = & 0  \\
y_{(1,7)}+10y_{(4,4)}+15y_{(5,3)}+6y_{(6,2)}  & = & 0  \\
y_{(2,6)}-6y_{(4,4)}-8y_{(5,3)}-3y_{(6,2)} & = & 0  \\
y_{(3,5)}+3y_{4,4)}+3y_{(5,3)}+y_{(6,2)} & =  & 0 \\
y_{(0,7)}+3y_{(1,6)}+3y_{(2,5)}+y_{(3,4)} & = & 0 \\
y_{(0,7)}+5y_{(1,6)}+10y_{(2,5)}+10y_{(3,4)}+5y_{(4,3)}+y_{(5,2)} & = &  0
\end{array}
\right. \fbox{$\begin{array}{cc}
x^1 & x^2  \\
x^1 & \bullet \\
x^1 & \bullet \\
x^1 & \bullet \\
\bullet & \bullet  \\
\bullet & \bullet 
\end{array}$}  \]
At order $8$, we have the parametric jets $\{ y_{(4,4)}, y_{(5,3)}, y_{(6,2)}, y_{(7,1)}, y_{(8,0)}\}$ with multiplicative variable $x^1$ {\it only} and thus ${\alpha}^2_8=1-1=0,{\alpha}^1_8=8-3=5$. At order 7, we have the 6 parametric jets $\{ y_{(2,5)}, y_{(3,4)}, y_{(4,3)}, y_{(5,2)}, y_{(6,1)}, y_{(7,0)} \}$ while at order $\leq 6$ we have $1+2+3+4+5+6+7=28$ parametric jets, that is a total of $28+6=34$ parametric jets with no multiplicative variable. Contrary to the Riquier/Janet/Gr\"{o}bner approach, we are sure that the only intrinsic number is the non-zero character ${\alpha}^1_8=5$ providing the parametric jets with only one multiplicative variable, in a coherent way with ([12],p 94) but in a completely different framework. As in the next example, the number of parametric jets with no multiplicative variable may however change a lot if we change the presentation of the module.   \\

\noindent
{\bf EXAMPLE 2.25}: With $n=4, m=1, q=2, k=\mathbb{Q}$, let us start with the unmixed perfect/radical polynomial ideal $\mathfrak{a}=({\chi}_1,{\chi}_2)\cap ({\chi}_3,{\chi}_4)\subset k[{\chi}_1,{\chi}_2, {\chi}_3, {\chi}_4]$ as in the end of [11] and consider the corresponding system made by $y_{13}=0, y_{14}=0, y_{23}=0, y_{24}=0$. This system is formally integrable because it is homogeneous but it is not evident to prove that it is also involutive. Integrating this system, it is easy to prove that the general solution $y=\varphi(x^1,x^2)+\psi (x^3,x^4)$ only depends on two arbitrary functions of two variables. However, {\it such a result has nothing to do with the CK theorem} because the quoted functions depend on different couples of variables. Hence it is not evident at first sight to know about the corresponding CK data. Let us make the following linear change of variables $x^1 \rightarrow x^1+x^4, x^2 \rightarrow x^2+x^3$ in order to obtain the involutive system:  \\
\[ \left\{  \begin{array}{lc}
y_{44}+y_{14} & = 0  \\
y_{34}+y_{13} & = 0  \\
y_{33}+y_{23}  & = 0  \\
y_{24}-y_{13}  & =  0
\end{array}
\right.  \fbox{$\begin{array}{cccc}
1 & 2 & 3 & 4  \\
1 & 2 & 3 & \bullet  \\
1 & 2 & 3 & \bullet  \\
1 & 2 & \bullet & \bullet
\end{array}$ }    \]
defining a $2$-pure differential module $M$.  \\
We have thus $par_2=\{y, y_1,y_2, y_3, y_4, y_{11}, y_{12}, y_{13}, y_{14}, y_{22}, y_{23}\}$ where $\{y,y_1, y_2, y_3, y_4\}$ have no multiplicative variable, $\{ y_{11}, y_{12}, y_{13}, y_{14} \}$ have the only multiplicative variable $x^1$ while $\{ y_{22}, y_{23}\}$ has the two multiplicative variables $(x^1, x^2)$. We have ${\alpha}^4_2=0, {\alpha}^3_2=0, {\alpha}^2_2=2, {\alpha}^1_2=4$ where only $\alpha={\alpha}^2_2$ and the last two series of two variables have an intrinsic meaning. \\
Setting $z^i=y_i$, we obtain the first order involutive system in reduced Spencer form:Ê\\
\[ \left\{  \begin{array}{l}
 d_4z^1-d_1z^4=0, d_4z^2-d_2z^4=0, d_4z^3-d_3z^4=0, d_4z^4+d_1z^4=0   \\
 d_3z^1-d_1z^3=0, d_3z^2-d_2z^3=0, d_3z^3+d_2z^3=0, d_3z^4+d_1z^3=0   \\
 d_2z^1-d_1z^2=0, d_2z^4-d_1z^3=0 
\end{array}
\right.  \fbox{$\begin{array}{cccc}
1 & 2 & 3 & 4  \\
1 & 2 & 3 & \bullet  \\
1 & 2 & \bullet & \bullet
\end{array}$ }    \]
with characters ${\alpha}^4_1=0, {\alpha}^3_1=0, \alpha={\alpha}^2_1=2(intrinsic), {\alpha}^1_1=4$ and the CK data for $z=g(x)$ are 
$\{ g^2(x^1,0,0,0), g^3(x^1,0,0,0), g^1(x^1,x^2,0,0), g^4(x^1,x^2,0,0) \}$ with now two series of two variables as before but also two series of one variable instead of the four found before. The system made by the equations of class $2$+class $3$ is involutive.  \\
Adding finally $y$ with $d_iy=z^i$, brings one additional equation to each class and does not therefore change the difference of characters, a result leading to the previous CK data to which one has to add $f(0,0,0,0)$. This is the standard way to get a differential module isomorphic to $M$ but defined by a first order system in Spencer form. A similar study can be done for the system $y_{44}=0, y_{34}=0, y_{33}=0, y_{24}-y_{13}=0$ coming from the primary ideal $\mathfrak{q}=(({\chi}_4)^2, {\chi}_3{\chi}_4, ({\chi}_3)^2, {\chi}_2{\chi}_4-{\chi}_1{\chi}_3)$ with radical $\mathfrak{p}=({\chi}_4,{\chi}_3)$ as in ([19], Ex 4.2). \\

\noindent
{\bf EXAMPLE 2.26}: In order to emphasize the importance of dealing with vector bundles in the differential geometric setting of this section and with differential fields or projective modules in the differential algebraic setting of the next section, we provide a tricky example of a linear system with coefficients in a true differential field which is not just a field of rational functions in the independent variables. With $n=2, m=1, q=2$, let us consider the non-linear system 
${\cal{R}}_2$:  \\
\[\left\{  \begin{array}{lc}
 y_{22}-\frac{1}{3}(y_{11})^3&=0\\
  y_{12}-\frac{1}{2}(y_{11})^2&=0  
  \end{array}
  \right.  \fbox{$\begin{array}{cc}
  1 & 2 \\
  1& \bullet
  \end{array}$}     \]
obtained by equating to zero two differential polynomials. Doing crossed derivatives, it is easy to check that the system is involutive and allows to define a true differential extension $K$ of $k=\mathbb{Q}$ which is isomorphic to $k(y,y_1,y_2,y_{11}, y_{111},... )$ if we set for example $d_2y_1=\frac{1}{2}(y_{11})^2$ and so on. By linearization, we get the following linear second order involutive system $R_2$ defined over $K$:  \\
\[\left\{  \begin{array}{lc}
 Y_{22}-(y_{11})^2Y_{11}&=0\\
  Y_{12}-y_{11}Y_{11}&=0  
  \end{array}
  \right.  \fbox{$\begin{array}{cc}
  1 & 2 \\
  1& \bullet
  \end{array}$}     \]
The various symbols of the first system are vector bundles over ${\cal{R}}_2$ while the symbols of the second system are vector spaces over $K$. As an exercise in order to understand the problems that may arise in general, we invite the reader to study similarly the non-linear second order system $y_{22}-\frac{1}{2}(y_{11})^2=0, y_{12}-y_{11}=0$ and conclude. The interested reader may look at ([14],VI.B.3 ,p 273 and VI.B.7 p 275) for criteria providing differential fields and based on the Spencer $2$-acyclicity property of the symbol at order $q$ ([14], Prop. III.1.3, p 92 and Theorem III.C.1, p 95).  \\

\noindent
{\bf 3  MACAULAY INVERSE SYSTEMS REVISITED} \\

Let $A$ be a {\it unitary ring}, that is $1,a,b\in A \Rightarrow a+b,ab \in A, 1a=a1=a$ and even an {\it integral domain} ($ab=0\Rightarrow a=0$ or $b=0$) with {\it field of fractions} $K=Q(A)$. However, we shall not always assume that $A$ is commutative, that is $ab$ may be different from $ba$ in general for $a,b\in A$. We say that $M={}_AM$ is a {\it left module} over $A$ if $x,y\in M\Rightarrow ax,x+y\in M, \forall a\in A$ or a {\it right module} $M_B$ over $B$ if the operation of $B$ on $M$ is $(x,b)\rightarrow xb, \forall b\in B$. If $M$ is a left module over $A$ and a right module over $B$ with $(ax)b=a(xb), \forall a\in A,\forall  b\in B, \forall x\in M$, then we shall say that $M={ }_AM_B$ is a {\it bimodule}. Of course, $A={ }_AA_A$ is a bimodule over itself. The category of left modules over $A$ will be denoted by $mod(A)$ while the category of right modules over $A$ will be denoted by $mod(A^{\it op})$. We define the {\it torsion submodule} $t(M)=\{x\in M\mid \exists 0\neq a\in A, ax=0\}\subseteq M$ and $M$ is a {\it torsion module} if $t(M)=M$ or a {\it torsion-free module} if $t(M)=0$. We denote by $hom_A(M,N)$ the set of morphisms $f:M\rightarrow N$ such that $f(ax)=af(x)$. In particular $hom_A(A,M)\simeq M$ because $f(a)=af(1)$ and we recall that a sequence of modules and maps is exact if the kernel of any map is equal to the image of the map preceding it. \\

When $A$ is commutative, $hom(M,N)$ is again an $A$-module for the law $(bf)(x)=f(bx)$ as we have $(bf)(ax)=f(bax)=f(abx)=af(bx)=a(bf)(x)$. In the non-commutative case, things are more complicate and we have:\\

\noindent
{\bf LEMMA 3.1}: Given ${}_AM$ and ${}_AN_B$, then $hom_A(M,N)$ becomes a right module over $B$ for the law $(fb)(x)=f(x)b$. Similarly, given 
${Ê}_AM_B$ and ${ }_AN$, then $hom_A(M,N)$ becomes a left module over $B$ for the law $(bf)(x)=f(xb)$.\\

\noindent
{\it Proof}:  In order to prove the first result we just need to check the two relations:
\[ (fb)(ax)=f(ax)b=af(x)b=a(fb)(x),\]
\[ ((fb')b")(x)=(fb')(x)b"=f(x)b'b"=(fb'b")(x).\]
The proof of the second result could be achieved similarly.  \\
\hspace*{12cm}               Q.E.D. \\
 
\noindent
{\bf DEFINITION 3.2}: A module $F$ is said to be {\it free} if it is isomorphic to a (finite) power of $A$ called the {\it rank} of $F$ over $A$ and denoted by 
$rk_A(F)$ while the rank of a module is the rank of a maximum free submodule. In the sequel we shall only consider {\it finitely presented} modules, 
namely {\it finitely generated} modules defined by exact sequences of the type $F_1 \stackrel{d_1}{\longrightarrow} F_0 \longrightarrow M\longrightarrow 0$ where $F_0$ and $F_1$ are free modules of finite ranks. For any short exact sequence $0\rightarrow M' \rightarrow M \rightarrow M" \rightarrow 0$, we have $rk_A(M)=rk_A(M')+rk_A(M")$. A module $P$ is called {\it projective} if there exists a free module $F$ and another (projective) module $Q$ such that $P\oplus Q\simeq F$. Accordingly, a {\it projective (free) resolution} of $M$ is a long exact sequence $... \stackrel{d_3}{\longrightarrow} P_2 \stackrel{d_2}{\longrightarrow} P_1 \stackrel{d_1}{\longrightarrow} P_0 \stackrel{p}{\longrightarrow} M \longrightarrow 0 $ where $P_0, P_1, P_2, ... $ are projective (free) modules, $M=coker(d_1)=P_0/ im(d_1)$ and $p$ is the canonical projection.  \\

A module $N$ over $A$ is {\it injective} if and only if $hom_A(\bullet,N)$ is an {\it exact functor}, that is transforms any short exact sequence into a short exact sequence or, equivalently ({\it Baer criterion}), if and only if any map $\mathfrak{a}\rightarrow N$ where $\mathfrak{a}\subset A$ is an ideal can be extended to a map $A \rightarrow N$. Accordingly, we may similarly define by duality an {\it injective resolution} of $M$ by using injective modules and reversing the arrows (See [24], p 67-74 for more details).  \\

Using the notation $M^*=hom_A(M,A)$, for any morphism $f:M\rightarrow N$, we shall denote by $f^*:N^*\rightarrow M^*$ the morphism which is defined by  $f^*(h)=h\circ f, \forall h\in hom_A(N,A)$ and satisfies $rk_A(f)=rk_A(im(f))=rk_A(f^*),\forall f\in hom_A(M,N)$(See [16], Corollary 5.3, p 179). We may take out $M$ in order to obtain the {\it deleted sequence} $... \stackrel{d_2}{\longrightarrow} P_1 \stackrel{d_1}{\longrightarrow} P_0 \longrightarrow 0$ and apply  $hom_A(\bullet,A)$ in order to get the sequence $... \stackrel{d^*_2}{\longleftarrow} P^*_1 \stackrel{d^*_1}{\longleftarrow} P^*_0 \longleftarrow 0$. \\

\noindent
{\bf DEFINITION 3.3}: A resolution of a short exact sequence $0 \rightarrow M' \stackrel{f}{\longrightarrow} M \stackrel{g}{\longrightarrow} M" \rightarrow 0 $ of $A$-modules is a short exact sequence $0 \rightarrow X' \stackrel{f}{\longrightarrow} X \stackrel{g}{\longrightarrow} X" \rightarrow 0 $ of exact complexes such that $X \stackrel{p}{\longrightarrow} M \rightarrow 0, X' \stackrel{p'}{\longrightarrow} M' \rightarrow 0, X" \stackrel{p"}{\longrightarrow} M" \rightarrow 0$ are resolutions and we shall say that the sequence of complexes is over the sequence of modules. Such a definition can also be used when the complexes are not exact and we have the long exact {\it connecting sequence} $... \rightarrow H_i(X) \rightarrow H_i(X") \rightarrow H_{i-1}(X') \rightarrow ... $ if we introduce the homology $H_i(X)$ of a decreasing complex $X_{i+1} \rightarrow X_i \rightarrow X_{i-1} $ with a similar result for the cohomology of increasing complexes. In particular, if any two are exact, the third is exact too ([15], Theorem II.1.15, p 196-203).  \\

\noindent
{\bf PROPOSITION 3.4}: The {\it extension modules}  $ext^0_A(M)=ker(d^*_1)=hom_A(M,A)=M^*$ and $ext^i_A(M)=ker(d^*_{i+1})/im(d^*_i), \forall i\geq 1$ do not depend on the resolution chosen and are torsion modules for $i\geq 1$. Using $hom_A(\bullet,N)$, one can similarly define $ext^i_A(M,N)$ with $ext^0_A(,N)=hom_A(M,N)$ and the $ext^i_A(M,N)$ vanish $\forall i>0$ whenever $M$ is a projective module or $N$ is an injective module (See [16] and [24] for the details). \\

Let $ A$ be a {\it differential ring}, that is a commutative ring with $n$ commuting {\it derivations} $\{{\partial}_1,...,{\partial}_n\}$ with ${\partial}_i{\partial}_j={\partial}_j{\partial}_i={\partial}_{ij}, \forall i,j=1,...,n$ such that ${\partial}_i(a+b)={\partial}_ia+{\partial}_ib$ and ${\partial}_i(ab)=({\partial}_ia)b+a{\partial}_ib, \forall a,b\in A$. More generally, a similar definition can be provided for a differential integral domain $A$ with unit $1\in A$ and will be used therafter whenever we shall need a {\it differential field} $\mathbb{Q}\subset K$ of coefficients that is a field ($a\in K\Rightarrow 1/a\in K$) with ${\partial}_i(1/a)=-(1/a^2){\partial}_ia$, for example in order to exhibit solved forms for systems of partial differential equations as in the preceding section. Using an implicit summation on multi-indices, we may introduce the (noncommutative) {\it ring of differential operators} $D=A[d_1,...,d_n]=A[d]$ with elements $P=a^{\mu}d_{\mu}$ such that $\mid \mu\mid<\infty$ and $d_ia=ad_i+{\partial}_ia$. The highest value of ${\mid}\mu {\mid}$ with $a^{\mu}\neq 0$ is called the {\it order} of the {\it operator} $P$ and the ring $D$ with multiplication $(P,Q)\longrightarrow P\circ Q=PQ$ is filtred by the order $q$ of the operators. We have the {\it filtration} $0\subset A=D_0\subset D_1\subset  ... \subset D_q \subset ... \subset D_{\infty}=D$. Moreover, it is clear that $D$, as an algebra, is generated by $A=D_0$ and $T=D_1/D_0$ with $D_1=A\oplus T$ if we identify an element $\xi={\xi}^id_i\in T$ with the vector field $\xi={\xi}^i(x){\partial}_i$ of differential geometry, but with ${\xi}^i\in A$ now. It follows that $D={ }_DD_D$ is a {\it bimodule} over itself, being at the same time a left $D$-module ${ }_DD$ by the composition $P \longrightarrow QP$ and a right $D$-module $D_D$ by the composition $P \longrightarrow PQ$ with $D_rD_s=D_{r+s}, \forall r,s \geq 0$ in any case. \\

If we introduce {\it differential indeterminates} $y=(y^1,...,y^m)$, we may extend $d_iy^k_{\mu}=y^k_{\mu+1_i}$ to ${\Phi}^{\tau}\equiv a^{\tau\mu}_ky^k_{\mu}\stackrel{d_i}{\longrightarrow} d_i{\Phi}^{\tau}\equiv a^{\tau\mu}_ky^k_{\mu+1_i}+{\partial}_ia^{\tau\mu}_ky^k_{\mu}$ for $\tau=1,...,p$. Therefore, setting $Dy^1+...+dy^m=Dy\simeq D^m$ and calling $I=D\Phi\subset Dy$ the {\it differential module of equations}, we obtain by residue the {\it differential module} or $D$-{\it module} $M=Dy/D\Phi$, denoting the residue of $y^k_{\mu}$ by ${\bar{y}}^k_{\mu}$ when there can be a confusion. Introducing the two free differential modules $F_0\simeq D^{m_0}, F_1\simeq D^{m_1}$, we obtain equivalently the {\it free presentation} $F_1\stackrel{d_1}{\longrightarrow} F_0 \rightarrow M \rightarrow 0$ of order $q$ when $m_0=m, m_1=p$ and $d_1={\cal{D}}=\Phi \circ j_q$. We shall moreover assume that ${\cal{D}}$ provides a {\it strict morphism} (see below) or, equivalently, that the corresponding system $R_q$ is formally integrable ([15]). It follows that $M$ can be endowed with a {\it quotient filtration} obtained from that of $D^m$ which is defined by the order of the jet coordinates $y_q$ in $D_qy$. We have therefore the {\it inductive limit} $0=M_{-1} \subseteq M_0 \subseteq M_1 \subseteq ... \subseteq M_q \subseteq ... \subseteq M_{\infty}=M$ with $d_iM_q\subseteq M_{q+1}$ but it is important to notice that $D_rD_q=D_{q+r} \Rightarrow D_rM_q= M_{q+r}, \forall q,r\geq 0 \Rightarrow M=DM_q, \forall q\geq 0$ {\it in this particular case}. It also follows from noetherian arguments and involution that 
$D_ rI_q=I_{q+r}, \forall r\geq 0$ though we have in general only $D_rI_s\subseteq I_{r+s}, \forall r\geq 0, \forall s<q$. It mut finally be noticed that the identification $(P_1,...,P_m)\leftrightarrow P_1y^1 + ... + P_my^m$ made by Piras in ([12], section 2, p 89) may bring the {\it rows} of the underlying differential operator of a system, in a coherent way with the identification $D^m\leftrightarrow Dy^1 + ... + Dy^m$ that we have used in the study of differential modules. As $A\subset D$, we may introduce the {\it forgetful functor} $for : mod(D) \rightarrow mod(A): { }_DM \rightarrow { }_AM$. In this paper, we shall go as far as possible with such an arbitrary differential ring $A$ though, in actual practice and thus in most of the examples considered, we shall use a differential field $K$ ([14]). We shall also assume that the ring $A$ is a noetherian ring (integral domain) in such a way that $D$ becomes a (both left and right) noetherian ring (integral domain).\\

More generally, introducing the successive CC as in the preceding section while changing slightly the numbering of the respective operators, we may finally obtain the {\it free resolution} of $M$, namely the exact sequence $\hspace{5mm} ... \stackrel{d_3}{\longrightarrow} F_2  \stackrel{d_2}{\longrightarrow} F_1 \stackrel{d_1}{\longrightarrow}F_0\longrightarrow M \longrightarrow 0 $. In actual practice, one must never forget that ${\cal{D}}=\Phi \circ j_q$ {\it acts on the left on column vectors in the operator case and on the right on row vectors in the module case}. Also, with a slight abuse of language, when ${\cal{D}}=\Phi \circ j_q$ is involutive as in section 2 and thus $R_q=ker( \Phi)$ is involutive, one should say that $M$ has an {\it involutive presentation} of order $q$ or that $M_q$ is {\it involutive} and $D_rM_q=M_{q+r}, \forall q,r \geq 0$ because $D_rD_q=D_{q+r}, \forall q,r \geq 0$. \\

In Section 2, the formal integrability of a system has been used in a crucial way in order to construct various differential sequences. Therefore, the algebraic counterpart provided by the next definition and proposition will also be used in a crucial way too in order to construct various resolutions of a differential module, though in a manner which is not so natural when dealing with applications to mathematical physics ([20], [21]). For this reason, we invite the reader to follow closely the arguments involved on the illustrating examples provided. To sart with, if $M$ and $N$ are two filtred differential modules and $f:M \rightarrow N$ is a differential morphism, that is a $D$-linear map with $f(Pm)=Pf(m), \forall P\in D$, then $f$ will be called an {\it homomorphism} of filtred modules if it induces $A$-linear maps $f_q=M_q \rightarrow N_q$. Chasing in the following commutative diagram:   \\
\[  \begin{array}{cccccl}
0  &   &  0  &  &  &     \\
\downarrow  &  & \downarrow  &  &  &  \\
M_q  & \stackrel{f_q}{\longrightarrow}  & N_q &  \rightarrow  &  coker(f_q) &  \rightarrow  0  \\
\downarrow &  \  & \downarrow  &  &  \downarrow &   \\
M  &  \stackrel{f}{\longrightarrow} & N & \rightarrow & coker(f) & \rightarrow  0
\end{array}   \]
while introducing $im(f)=I\subseteq N, im(f_q)=I_q\subseteq N_q$, we may state: \\

\noindent
{\bf DEFINITION 3.5}: A differential morphism $f$ is said to be a {\it strict homomorphism} if the two following equivalent properties hold:  \\
1) There is an induced monomorphism $0\rightarrow coker(f_q) \rightarrow coker(f), \forall q\geq 0$.  \\
2) $f_q(M_q)=f(M)\cap N_q $, that is $ I_q=I\cap N_q$.\\
A sequence made by strict morphisms will be called a {\it strict sequence}. In order to fulfill the conditions of the definition, it is most of the time necessary to "{\it shift} " the filtration of a differential module $M$ by setting $M(r)_q=M_{q+r}$ in such a way that $q$ could be negative and we shall therefore always assume that $M_q=0, \forall q\ll 0$. \\

\noindent
{\bf PROPOSITION 3.6}: If we have a strict short exact sequence $0 \rightarrow M' \stackrel{f}{\longrightarrow} M \stackrel{g}{\longrightarrow} M" \rightarrow 0$ in which $\exists q\gg 0$ such that $D_rM_q=M_{q+r}, \forall r\geq 0$, then $\exists q',q" \gg0$ such that $D_rM'_{q'}=M'_{q'+r}, D_rM"_{q"}=M"_{q"+r}, \forall r\geq 0$ and conversely. We may thus assume that $q=q'=q"$ in both cases by choosing $q \gg 0$. More generally, an exact sequence of filtred differential modules is strictly exact if and only if the associated sequence of graded modules is exact in a way dualizing the differential geometric framework, on the condition to shift conveniently the various filtrations involved.\\

\noindent
{\it Proof}: First of all, setting $G=gr(M), G'=gr(M'), G"=gr(M")$, we have the commutative and exact diagram:  \\
\[  \begin{array}{rcccccl}
  &  0  & &  0  & &   0  &   \\
  & \downarrow  &  & \downarrow  &  & \downarrow &    \\
0 \rightarrow & M'_{q-1} & \stackrel{f_{q-1}}{\longrightarrow}  & M_{q-1} & \stackrel{g_{q-1}}{\longrightarrow} & M"_{q-1} & \rightarrow 0 \\
  & \downarrow  &  & \downarrow  &  & \downarrow &    \\
0 \rightarrow & M'_q & \stackrel{f_q}{\longrightarrow}  & M_q & \stackrel{g_q}{\longrightarrow} & M"_q & \rightarrow 0 \\
& \downarrow  &  & \downarrow  &  & \downarrow &    \\
0 \rightarrow & G'_q & \stackrel{gr_q(f)}{\longrightarrow}  & G_q & \stackrel{gr_q(g)}{\longrightarrow} & G"_q & \rightarrow 0 \\
& \downarrow  &  & \downarrow  &  & \downarrow &    \\
 & 0 &  & 0  & &  0 & 
\end{array}  \]
Indeed, as $g$ is a strict epimorphism, it follows that $g_q$ is surjective $\forall q\geq 0$. Also, as $f$ is a monomorphism, then $f_q$ is also a monomorphism $\forall q \geq 0$ by restriction. Moreover, as $f$ is also strict, we obtain successively by chasing:  \\
\[     ker(g_q)=f(M')\cap M_q=f(M'_q)=f_q(M'_q)=im(f_q)  \]
It follows that the two upper rows are exact and the bottom row is thus exact too $\forall q\geq 0$ from the snake theorem in homological algebra ([2],[15],[16], [25]).  \\
This result provides the short exact sequence $0 \rightarrow G' \stackrel{gr(f)}{\longrightarrow} G \stackrel{gr(g)}{\longrightarrow} G" \rightarrow 0 $ of graded modules. \\
Let us now consider the following commutative diagram with maps such as $\xi \otimes m \rightarrow \xi m$ and where the upper row is exact because $D_1\simeq A \oplus T$ is free over $A$:\\
\[  \begin{array}{rcccccl}
0 \rightarrow & D_1{\otimes}_AM'_q & \stackrel{f}{\longrightarrow} & D_1{\otimes}_A M_q  & \stackrel{g}{\longrightarrow} & D_1{\otimes}_AM"_q & \rightarrow 0  \\
  & \downarrow & & \downarrow & & \downarrow &  \\
  0 \rightarrow & M'_{q+1} & \stackrel{f_{q+1}}{\longrightarrow} & M_{q+1} & \stackrel{g_{q+1}}{\longrightarrow} & M"_{q+1}  & \rightarrow 0  
  \end{array}  \]

If the central map is surjective, then the map on the right is also surjective, that is $D_1M_q=M_{q+1} \Rightarrow D_1M"_q=M"_{q+1}$ and thus $q=q"$. This is the typical situation met in a finite presentation of a system already considered. Moreover, $D_1M'_q\subseteq M'_{q+1} \Rightarrow TG'_q\subseteq G'_{q+1}$ like in the following commutative and exact diagrams where the left one is holding for a (formally integrable) system while the corresponding right one is holding for an arbitrary filtred module $M$ with $gr(M)=G$:  \\
\[ \begin{array}{rccclcrcccl}
 &   0  &   &  0  &                                               \hspace{1cm}       &                     &     &   0  &   &  0  &                                              \\                                         
    &  \downarrow  &  &  \downarrow  &                                           &     \hspace{1cm} & & \uparrow  &  &  \uparrow  &                    \\
0 \rightarrow & g_{q+1}  &  \rightarrow  &  T^* \otimes R_q &      &     \hspace{1cm} &     & T{\otimes}_AM_q  &  \rightarrow  &  G_{q+1} &  \rightarrow 0 \\
    &   \downarrow  &  &  \downarrow  &                                        &  \hspace{1cm}&   &  \uparrow  &  &  \uparrow  &                \\
 0 \rightarrow &R_{q+1} & \rightarrow & J_1(R_q) &               &   \hspace{1cm} &    &  D_1{\otimes}_AM_q & \rightarrow & M_{q+1} &    \\
  &  \downarrow  &  &  \downarrow  &                                             & \hspace{1cm}  &   &\uparrow  &  &  \uparrow  &               \\
0 \rightarrow &    R_q  &  =  & R_q  & \rightarrow 0                     &  \hspace{1cm}  &   0 \rightarrow &    M_q  &  =  & M_q  & \rightarrow 0    \\
   &  \downarrow  &  &  \downarrow  &                                           & \hspace{1cm} & &   \uparrow  &  &  \uparrow  &       \\
     &   0  &  &  0  &                                                                              &  \hspace{1cm} &  &   0  &  &  0  &    
      \end{array}       \]
 In these diagrams, the upper morphism is the composition $g_{q+1} \stackrel{\delta}{\longrightarrow} T^*\otimes g_q \rightarrow T^*\otimes R_q$ in the system diagram and the composition $T{\otimes}_AM_q \rightarrow T{\otimes}_AG_q \rightarrow G_{q+1}$ in the module diagram. Accordingly, a chase is showing that $D_1M_q=TM_q+M_q\subseteq M_{q+1}$ with equality if and only if  $TG_q=G_{q+1}$.\\
From noetherian arguments for polynomial rings in commutative algebra, it follows that $G'$ is finitely generated and we may choose for $q'$ the maximum order of a minimum set of generators.\\
 Conversely, if $D_rM'_{q'}=M'_{q'+r}, D_rM"_{q"}=M"_{q"+r}, \forall r\geq 0$, we may choose $q=sup(q',q")$ and we have thus $D_1M'_q=M'_{q+1}, D_1M"_q=M"_{q+1}\Rightarrow D_1M_q=M_{q+1}$, using again the snake theorem. \\
 As a byproduct, it is always possible to find $q\gg 0$ such that we could have {\it at the same time} $D_rM_q=M_{q+r}, D_rM'_q=M'_{q+r}, D_rM"_q=M"_{q+r}, \forall r\geq 0$ in the two situations considered.  \\
 We end this proof with a comment on the prolongation of symbols and graded modules which, in our opinion based on more than thirty years spent on computing and applying these dual concepts, is not easy to grasp. For this, let us consider the corresponding diagrams:  \\
 \[  \begin{array}{rccccccccl}
   & 0  &  &  0  &                                                                                     \hspace{1cm}                          &  0 &  & 0 &   \\
  & \downarrow &  & \downarrow &                                                               \hspace{1cm}       & \uparrow   &       & \uparrow  &    \\
0 \rightarrow & g_{q+1} & \rightarrow &  S_{q+1}T^* \otimes E  &\hspace{1cm}    &S_{q+1}T{\otimes}_AE^* & \rightarrow & G_{q+1} & \rightarrow 0  \\
   &\hspace{3mm}\downarrow\delta &  &\hspace{3mm}\downarrow\delta  &   \hspace{1cm} &\hspace{5mm} \uparrow{\delta}^*  &  &\hspace{5mm}  \uparrow {\delta}^* &   \\
0\rightarrow & T^*\otimes g_q & \rightarrow & T^*\otimes S_qT^*\otimes E &\hspace{1cm}     &T{\otimes}_AS_qT{\otimes}_A E^*  & \rightarrow & T{\otimes}_AG_q&\rightarrow 0 
\end{array}   \] 
Indeed, exactly as we have in general $R_{q+1}\subseteq {\rho}_1(R_q)\Rightarrow g_{q+1}\subseteq {\rho}_1(g_q)$, there is no corresponding concept in module theory without a reference to a presentation. In the differential geometric framework, ${\rho}_1(g_q)$ is the {\it reciprocal image} of $\delta$, that is the subset (not always a vector bundle !) of $S_{q+1}T^*\otimes E$ made by elements having an image in $T^*\otimes g_q$ under $\delta$. Such a definition is also the one of "{\it fiber product}" in ([8], III.5, p 88-91). \\
 \hspace*{12cm}       Q.E.D.    \\

\noindent
{\bf EXAMPLE 3.7}: Though this is not evident at first sight when $m=1, n=2, A=\mathbb{Q}[x^1,x^2]$, we invite the reader to prove that the third order linear system $y_{222}+x^2y_2=0, y_{111}+y_2-y=0$ has the same formal solutions as the third order system $y_{111}-y=0, y_2=0$ which is defined over $\mathbb{Q}$, a result leading to the generating involutive third order linear system $y_{222}=0, y_{122}=0, y_{112}=0, y_{111}-y=0, y_{22}=0, y_{12}=0,y_2=0$. We have $M_0=\{\bar{y}\},  M_1=\{\bar{y},{\bar{y}}_1\}, M_2=\{\bar{y},{\bar{y}}_1,{\bar{y}}_{11}\}=M$ while using only parametric jets because $d_1{\bar{y}}_{11}={\bar{y}}_{111}=\bar{y}$ and thus $D_1I_1=I_2, D_1I_2\subset I_3$ with a strict inclusion, $DI_3=I$. (See the similar Examples III.2.64 and III.3.11 in [15] for the details). \\

Roughly speaking, homological algebra has been created for finding intrinsic properties of modules not depending on their presentations or even on their resolutions and we now exhibit another approach by defining the {\it formal adjoint} of an operator $P$ and an operator matrix $\cal{D}$:  \\

\noindent
{\bf DEFINITION 3.8}: Setting $P=a^{\mu}d_{\mu}\in D  \stackrel{ad}{\longleftrightarrow} ad(P)=(-1)^{\mid\mu\mid}d_{\mu}a^{\mu}   \in D $, we have $ad(ad(P))=P$ and $ad(PQ)=ad(Q)ad(P), \forall P,Q\in D$. Such a definition can be extended to any matrix of operators by using the transposed matrix of adjoint operators and we get:  
\[ <\lambda,{\cal{D}} \xi>=<ad({\cal{D}})\lambda,\xi>+\hspace{1mm} {div}\hspace{1mm} ( ... )  \]
from integration by part, where $\lambda$ is a row vector of test functions and $<  > $ the usual contraction.  \\

\noindent
{\bf LEMMA 3.9}: If $f\in aut(X)$ is a local diffeomorphisms on $X$, we may set $ x=f^{-1}(y)=g(y)$ and we have the {\it identity}:
\[   \frac{\partial}{\partial y^k}(\frac{1}{\Delta (g(y))} {\partial}_if^k(g(y))\equiv 0.   \]

\noindent
{\bf PROPOSITION 3.10}: If we have an operator $E\stackrel{\cal{D}}{\longrightarrow} F$, we may obtain by duality an operator ${\wedge}^nT^*\otimes E^*\stackrel{ad(\cal{D})}{\longleftarrow} {\wedge}^nT^*\otimes F^*$. \\

\noindent
{\bf EXAMPLE 3.11}: In order to understand how the Lemma is involved in the Proposition, let us revisit relativistic electromagnetism (EM) in the light of these results when $n=4$. First of all,  we have $dA=F \Rightarrow dF=0$ in the sequence ${\wedge}^1T^*\stackrel{d}{\longrightarrow} {\wedge}^2T^* \stackrel{d}{\longrightarrow} {\wedge}^3T^*$ and the {\it field equations} of EM (first set of Maxwell equations) are invariant under any local diffeomorphism $f\in aut(X)$. By duality, we get the sequence ${\wedge}^4T^*\otimes {\wedge}^1T \stackrel{ad(d)}{\longleftarrow} {\wedge}^4T^*\otimes {\wedge}^2T \stackrel{ad(d)}{\longleftarrow} {\wedge}^4T^*\otimes {\wedge}^3T$ which is locally isomorphic (up to sign) to ${\wedge}^3T^* \stackrel{d}{\longleftarrow} {\wedge}^2T^* \stackrel{d}{\longleftarrow} {\wedge}^1T^*$ and the {\it induction equations} ${\partial}_i{\cal{F}}^{ij}={\cal{J}}^j$ of EM (second set of Maxwell equations) are thus also invariant under any $f\in aut(X)$. Indeed, using the last lemma and the {\it identity} ${\partial}_{ij}f^l{\cal{F}}^{ij}\equiv 0$, we have: \\
\[\frac{\partial}{\partial y^k}(\frac{1}{\Delta}{\partial}_i f^k{\partial}_j f^l{\cal{F}}^{ij})=\frac{1}{\Delta} {\partial}_i f^k \frac{\partial}{\partial y^k}({\partial}_j f^l{\cal{F}}^{ij})=\frac{1}{\Delta}{\partial}_i({\partial}_j f^l{\cal{F}}^{ij})=\frac{1}{\Delta}{\partial}_j f^l{\partial}_i{\cal{F}}^{ij} \]
Accordingly, it is not correct to say that the conformal group is the biggest group of invariance of Maxwell equations in physics as it is only the biggest group of invariance of the Minkowski constitutive laws in vacuum [21]. Finally, according to Proposition 3.3 both sets of equations can be parametrized {\it independently}, the first by the potential, the second by the so-called pseudopotential (See [15], p 492 for more details).\\

Now, with operational notations, let us consider the two differential sequences:  \\
\[   \xi  \stackrel{{\cal{D}}}{\longrightarrow} \eta \stackrel{{\cal{D}}_1}{\longrightarrow} \zeta  \]
\[   \nu  \stackrel{ad({\cal{D}})}{\longleftarrow} \mu \stackrel{ad({\cal{D}}_1)}{\longleftarrow} \lambda   \]
where ${\cal{D}}_1$ generates all the CC of ${\cal{D}}$. Then ${\cal{D}}_1\circ {\cal{D}}\equiv 0 \Longleftrightarrow ad({\cal{D}})\circ ad({\cal{D}}_1)\equiv 0 $ but $ad({\cal{D}})$ may not generate all the CC of $ad({\cal{D}}_1)$. Passing to the module framework, we just recognize the definition of $ext^1_D(M)$. Now, exactly like we defined the differential module $M$ from $\cal{D}$, let us define the differential module $N$ from $ad(\cal{D})$. Then $ext^1_D(N)=t(M)$ does not depend on the presentation of $M$. More generally, changing the presentation of $M$ may change $N$ to $N'$ but we have ([15],p 651)([16],p 203):  \\

\noindent
{\bf THEOREM 3.12}: The modules $N$ and $N'$ are {\it projectively equivalent}, that is one can find two projective modules $P$ and $P'$ such that $N\oplus P\simeq N' \oplus P'$ and we obtain therefore $ext^i_D(N)\simeq ext^i_D(N'), \forall i\geq 1$.  \\

Having in mind that $D$ is a $A$-algebra, that $A$ is a left $D$-module with the standard action $(D,A) \longrightarrow A:(P,a) \longrightarrow P(a):(d_i,a)\longrightarrow {\partial}_ia$ and that $D$ is a bimodule over itself, {\it we have only two possible constructions leading to the following two definitions}:  \\

\noindent
{\bf DEFINITION 3.13}: We define the {\it inverse system} $R=hom_A(M,A)$ of $M$ and set $R_q=hom_A(M_q,A)$ as the {\it inverse system of order} $q$. From the injective limit of the filtration of $M$ we deduce the {\it projective limit} $R=R_{\infty} \longrightarrow ... \longrightarrow R_q \longrightarrow ... \longrightarrow R_1 \longrightarrow R_0$. It follows that $f_q\in R_q:y^k_{\mu} \longrightarrow f^k_{\mu}\in A$ with $a^{\tau\mu}_kf^k_{\mu}=0$ defines a {\it section at order} $q$ and we may set $f_{\infty}=f\in R$ for a {\it section} of $R$. For a ground field of constants $k$, this definition has of course to do with the concept of a formal power series solution. However, for an arbitrary differential ring $A$ or differential field $K$, {\it the main novelty of this new approach is that such a definition has nothing to do with the concept of a formal power series solution} as illustrated in the next examples. Nevertheless, in actual practice, it is always simpler to deal with a differential field $K$ in order to have finite dimensional vector spaces at each order $q$ for applications.\\

Now, if $A,B$ are rings and ${ }_AM, { }_BL_A, { }_BN$ are modules, using the second part of Lemma 3.1 and the relation $l\otimes am=la\otimes m$ with the left action $b(l\otimes m)=bl\otimes m,\forall a\in A, \forall b\in B,\forall l\in L, \forall m\in M$, we may provide the so-called {\it adjoint isomorphism} as in ([25], Th 2.11, p. 37):  \\ 
\[   \varphi :hom_B(L{\otimes}_A M, N)\stackrel{\simeq}{\longrightarrow} hom_A(M, hom_B(L,N))   \] 
saying that there is a one-to-one correspondence between maps of the form $L\otimes M \rightarrow N$ and maps of the form $M \rightarrow hom(L,N)$ or, fixing an element $m\in M$, providing a parametrization set of maps of the form $L \rightarrow N$ in both cases.  \\
With more details, starting with a $B$-morphism $f:L{\otimes}_AM \rightarrow N$, one may define a $A$-morphism $\varphi(f):M \rightarrow hom_B(L,N)$ by the formula $(\varphi(f)(m))(l)=f(l\otimes m)$. It follows that such a $\varphi$ is a monomorphism because it is defined on the basis of simple tensors in $L{\otimes}_AM$ and it remains to check that it is also an epimorphism by constructing an inverse $\psi$. For this, starting with a $A$-morphism $g:M\rightarrow hom_B(L,N)$, we just define $\psi(g)=f$ by $f(l\otimes m)=(g(m))(l)$. We have in particular: \\
\[  (\varphi(f)(am))(l)=f(l\otimes am)=f(la\otimes m)=(\varphi(f)(m))(la)=a(\varphi(f)(m))(l)  \]
and thus $\varphi(f)(am)=a(\varphi(f)(m))$ in a coherent way with $hom_A$ and Lemma 3.1. \\

\noindent
With $M={ }_DM,L={ }_AD_D,N={ }_AA$ and $1\in A\subset D$, one obtains the isomorphism:  \\
\[  hom_A(M,A)=hom_A(D{\otimes}_  DM,A)\simeq hom_D(M,hom_A(D,A))\]
where $(\varphi(f)(m))(1)=f(m)$ if we identify $m$ with $1\otimes m$ and $D^*=hom_A(D,A)$ is an injective module because of Baer's criterion when $A=K=k$ is a field of constants, that is when $D$ is a commutative ring ([25], Th 3.20, p. 67) or ([2], Proposition 11, p 18)(See [24], in particular chapter 4,  for more details). It follows that $ext^i_D(M,D^*)=0, \forall i>0$ (See section 2.4 in [7]). However, though $M$ is a left $D$-module by assumption and $hom_A(D,A)$ is also a left $D$-modules with $(Qh)(P)=h(PQ), \forall h\in hom_A(D,A),\forall P,Q\in D$ because of Lemma 3.1, there is no similar reason " {\it a priori} " that $hom_A(M,A)$ could {\it also} be a left $D$-module. Moreover, we notice that $(ah)(P)=h(Pa)\neq h(aP)=a(h(P)), \forall a\in A$ unless $A$ is a field $k$ of constants.\\

\noindent
ACCORDINGLY, THIS APPROACH IS NOT CONVENIENT AND MUST BE MODIFIED WHEN $A$ IS A TRUE DIFFERENTIAL RING OR $K$ IS A TRUE DIFFERENTIAL FIELD, THAT IS WHEN $D$ IS NOT COMMUTATIVE.  \\

\noindent
{\bf DEFINITION 3.14}: We may define the right differential module $M^*=hom_D(M,D)$.  \\

The next crucial theorem will allow to provide the module counterpart of the differential geometric construction of the Spencer operator provided in Section 2 (Compare to [2] and [15]). For a more general approach, we shall consider a differential ring $A$ with unity $1$ and set $D=A[d]$ as in ([15], chapter IV). \\

\noindent
{\bf THEOREM 3.15}: When $M$ and $N$ are left $D$-modules, then $hom_A(M,N)$ and $M{\otimes}_AN$ are left $D$-modules. In particular $R=hom_A(M,A)$ is also a left $D$-module for the Spencer operator. Moreover, the structures of left $D$-modules existing therefore on $M{\otimes}_AN$ and $hom_A(N,L)$ are now coherent with the {\it adjoint isomorphism} for $mod(D)$:  \\
\[   \varphi :  hom_D(M{\otimes}_AN,L) \stackrel{\simeq}{\longrightarrow} hom_D(M,hom_A(N,L)) \hspace{5mm} ,\forall L,M,N\in mod(D)    \]
Finally, if $M$ and $N$ are right $D$-modules, then $hom_A(M,N)$ is a left $D$-module. Moreover, if $M$ is a left $D$-module and $N$ is a right $D$-module, then $M{\otimes}_AN$ is a right $D$-module. It follows that we have also $R=hom_A(M,A)\simeq hom_D(M,D^*)$ but in a quite different framework.\\

\noindent
{\it Proof}:  Let us define for any $f\in hom_A(M,N)$:   \\
\[   (af)(m)=af(m)=f(am) \hspace{1cm} \forall a\in A, \forall m\in M\]
\[   (\xi f)(m)=\xi f(m)-f(\xi m)  \hspace{1cm}  \forall \xi ={\xi}^id_i\in T, \forall m\in M  \]
It is easy to check that $\xi a=a \xi+\xi (a)$ in the operator sense and that $\xi\eta -\eta\xi =[\xi,\eta]$ is the standard bracket of vector fields. We have in particular with $d$ in place of any $d_i$: \\
\[  \begin{array}{rcl}
((da)f)(m)=(d(af))(m)=d(af(m))-af(dm)&=&(\partial a)f(m)+ad(f(m))-af(dm)\\
       &=& (a(df))(m)+(\partial a)f(m)  \\
       &=& ((ad+\partial a)f)(m)
       \end{array}  \]
 We may then define for any $m\otimes n\in M{\otimes}_AN$ with arbitrary $m\in M$ and $n\in N$:   \\
 \[      a(m\otimes n)=am\otimes n=m\otimes an\in M{\otimes}_AN  \]
 \[  \xi (m\otimes n)=\xi m\otimes n + m\otimes \xi n \in M{\otimes}_AN   \]
 and conclude similarly with:   \\
 \[  \begin{array}{rcl}
  (da)(m\otimes n)=d(a(m\otimes n)) & = & d(am\otimes n)\\
                         & = &  d(am)\otimes n+am\otimes dn  \\
                             & = & (\partial a)m\otimes n + a(dm)\otimes n + am\otimes dn  \\
                                & = & (ad+\partial a)(m\otimes n)
                                \end{array}    \]
Using $A$ or $K=Q(A)$ in place of $N$, we finally get $(d_if)^k_{\mu}=(d_if)(y^k_{\mu})={\partial}_if^k_{\mu}-f^k_{\mu +1_i}$ that is {\it we recognize exactly the Spencer operator} that we have used in the second part and thus:\\
\[  (d_i(d_jf))^k_{\mu}={\partial}_{ij}f^k_{\mu}-({\partial}_if^k_{\mu+1_j}+{\partial}_jf^k_{\mu+1_i})+f^k_{\mu+1_i+1_j} \Rightarrow d_i(d_jf)=d_j(d_if)=d_{ij}f \]
In fact, $R$ is the {\it projective limit} of ${\pi}^{q+r}_q:R_{q+r}\rightarrow R_q$ in a coherent way with jet theory ([18],[19]). In the more specific case of $D^*=hom_A(D,A)$, the upper index $k$ is not present and we have thus $(af)_{\mu}=af_{\mu}$ with $ (d_if)_{\mu}={\partial}_if_{\mu}-f_{\mu+1_i},\forall f\in D^*,\forall a\in A,\forall i=1,...,n$ (Compare to [23], chapter 4, where the Spencer operator is lacking). This left $D$-module structure on $D^*$ is quite different from the one provided by Lemma 3.1 but coincide with it up to sign when $A=k$.\\
With $\varphi(f)=g$, the third result is entrelacing the two left structures that we have just provided through the formula $(g(m))(n)=f(m\otimes n)\in N$ defining the map $\varphi$ whenever $f\in hom_D(M{\otimes}_A N,L)$ is given. Using any $\xi\in T$, we get successively in $L$ (Compare to [2], Proposition 2.1.3, p 54):  \\
\[  \begin{array}{rcl}
(\xi(g(m)))(n)& = & \xi((g(m))(n))-(g(m))(\xi n)  \\
                                   & = & \xi(f(m\otimes n))-f(m\otimes \xi n)    \\
                                   & = & f(\xi(m\otimes n))-f(m\otimes \xi n)  \\
                                   & = & f(\xi m\otimes n+m\otimes \xi n)-f(m\otimes \xi n)  \\
                                   & = & f(\xi m\otimes n)  \\
                                   & = & (g(\xi m))(n)
\end{array}  \]
and thus $ \xi(g(m))=g(\xi m), \forall m\in M $ or simply $\xi\circ g=g\circ \xi$.   \\
For any $g\in hom_D(M,hom_K(N,L))$, we may define the inverse $\psi$ of $\varphi$ through the formula $\psi(g)(m\otimes n)=(g(m))(n)\in L$ by checking the bilinearity over $A$ of $(m,n) \rightarrow (g(m))(n)$ and studying as before the action of any $\xi\in T$.\\
Finally, if $M$ and $N$ are right $D$-modules, we just need to set $(\xi f)(m)=f(m\xi)-f(m)\xi, \forall \xi\in T, \forall m\in M$ and conclude as before. Similarly, if $M$ is a left $D$-module and $N$ is a right $D$-module, we just need to set $(m\otimes n)\xi=m\otimes n\xi - \xi m \otimes n$. \\
The last result is even more tricky and we provide two different proofs. \\
If $M$ is finitely presented, applying $hom_D(\bullet,D^*)$ to a free presentation $D^p\stackrel{{\cal{D}}}{\longrightarrow}D^m \rightarrow M \rightarrow 0 $, we obtain the exact sequence $0\rightarrow hom_D(M,D^*) \rightarrow D^{*m} \rightarrow D^{*p}$ because $hom_D(D,D^*)=D^*$. As any module over $D$ is a module over $A$, applying $hom_A(\bullet, A)$ to the same sequence, we get the exact sequence $0 \rightarrow hom_A(M,A) \rightarrow D^{*m} \rightarrow D^{* p}$ and thus an isomorphism $R=hom_A(M,A)\simeq hom_D(M, D^*)$.\\
More generally, when $M={  }_DM,N={ }_DD_D$ and $L={ }_DA$, it follows that we have an isomorphism $hom_D(M{\otimes}_AD,A)\simeq hom_D(M,D^*)$. Let us construct an isomorphism $M{\otimes}_AD\simeq D{\otimes}_AM$ of bimodules over $D$. Indeed, we may use for $M{\otimes}_AD$ the left structure over $D$ provided by the previous result with $M={ }_DM$ and $D={ }_DD$ while introduing a right structure over $D$ by defining $(m\otimes P)Q=m\otimes PQ$. Similarly, we may use for $D{\otimes}_AM$ the right structure provided by the last result with $D=D_D$ and $M={ }_DM$ while introducing a left structure by defining $Q(P\otimes m)=QP\otimes m$. We obtain for example: \\
\[ Q((P\otimes m)\eta)=Q(P\eta\otimes m-P\otimes \eta m)=QP\eta\otimes m-QP\otimes \eta m=(Q(P\otimes m))\eta  \]
The isomorphism is obtained by setting (Compare to [1], Prop. 2.2.8 and [26], Prop. 4.1.3): \\
\[    P\otimes m=P(1\otimes m) \longrightarrow  P(1\otimes m)  \]
\[  m\otimes P=(m\otimes 1)P  \longrightarrow  (1\otimes m)P   \]
and checking that :  \\
\[ Q(P\otimes m)=QP\otimes m \longrightarrow QP(m\otimes 1)=Q(P(m\otimes 1))  \]
\[  (m\otimes P)Q==(m\otimes 1)PQ\longrightarrow (1\otimes m)PQ=((1\otimes m)P)Q  \]
both with:  \\
\[ \begin{array}{lcl}
 \xi (m\otimes P)=(\xi m \otimes 1)P+(m \otimes 1) \xi P &  \rightarrow & (1 \otimes \xi m)P+(1 \otimes m)\xi P  \\
                                                                                                   &  = & (1 \otimes \xi m)P +(\xi \otimes m - 1 \otimes \xi m)P  \\
                                                                                                   &  = & (\xi \otimes m )P
 \end{array}     \]
and thus $\xi(m \otimes 1)P \rightarrow \xi (1 \otimes m)P $, that is $m\otimes 1 \rightarrow 1\otimes m$. \\
Using an induction on $ord(P)$, we have successively: \\
\[ \begin{array}{rcl}
P\xi \otimes m= P\xi (1\otimes m) &\rightarrow &P(\xi (m \otimes 1))=P(\xi m \otimes 1+ m \otimes \xi)=P(\xi m \otimes 1) + P(m \otimes 1) \xi\\
   & \rightarrow & P(1 \otimes \xi m) + P((\xi \otimes m) - (1 \otimes \xi m))= P(\xi \otimes m)  
   \end{array}    \]
and obtain therefore a functorial isomorphism $D {\otimes}_A M \simeq M {\otimes}_AD$ of bimodules. \\
It follows that $hom_D(M,hom_A(D,A))\simeq hom_D(M{\otimes}_AD,A)\simeq hom_D(D{\otimes}_AM,A) \simeq hom_A(M,A)$ because $a \otimes m=1 \otimes am, \forall a\in A \subset D, \forall m \in M$.  \\
With more details, if we set $\tilde{\varphi}(\tilde{f})=g $ and define $\tilde{\varphi}:hom_D(D{\otimes}_AM,A)\longrightarrow hom_D(M,D^*)$ by $(g(m))(P)=\tilde{f}((1\otimes m)P)$ when $\tilde{f}\in hom_D(D{\otimes}_AM,A)$, we get successively:  \\
\[   \begin{array}{lcl}
(\xi (g(m)))(P) & = & \xi ((g(m))(P))-(g(m))(\xi P)  \\
                                 & = & \xi \tilde{f}((1\otimes m)P)-\tilde{f}((1\otimes m)\xi P)   \\
                                 & = &  \tilde{f}(\xi (1\otimes m)P)-\tilde{f}((1\otimes m)\xi P)  \\
                                 & = & \tilde{f}(\xi (1 \otimes m)P)-\tilde{f}((\xi \otimes m)P)-(1 \otimes \xi m)P) \\
                                 & = & \tilde{f}((\xi\otimes m)P)-\tilde{f}((\xi \otimes m)P)+\tilde{f}((1\otimes \xi m)P)  \\
                                 & = & \tilde{f}((1\otimes \xi m )P)  \\
                                 & = & (g(\xi m))(P)
 \end{array}    \]
 and thus  $\xi (g(m))=g(\xi m)$ or simply $\xi \circ g= g\circ \xi$ again.  \\
 It is finally important to notice that the left and right $D$-structures that can be given to ${ }_DD{\otimes}_A{ }_DM$ and $D_D{\otimes}_A{ }_DM$ respectively do not allow to provide a bimodule structure on $D{\otimes}_AM$. Indeed, we have on one side: \\
 \[  \begin{array}{lcl}
 (\xi(P\otimes m))\eta &=& (\xi P\otimes m+P \otimes \xi m)\eta  \\
                                      & = & \xi P \eta \otimes m-\xi P\otimes \eta m + P \eta\otimes \xi m -P \otimes \eta \xi m  
                                      \end{array}   \] 
while we have on the other side:   \\
\[   \begin{array}{lcl}
\xi((P\otimes m)\eta)  & = & \xi (P\eta \otimes m-P\otimes \eta m)  \\
                                      & = & \xi P \eta \otimes m +P \eta \otimes \xi m -\xi P \otimes \eta m -P \otimes \xi \eta m
                                      \end{array}  \]
and thus $(\xi(P \otimes m))\eta - \xi((P\otimes m)\eta )=P\otimes [\xi,\eta]m\neq 0$ unless $D$ is a commutative ring.\\
                                 
\hspace*{12cm}     Q.E.D.   \\

When $A$ is an integral domain with field of fractions $K=Q(A)$, then it is known from ([2], See Example 1, p 18) that $K$ is an injective $A$-module. It thus follows from the previous results and the diagram of ([2], Proposition 11, p 18) that $hom_A(D,K)$ is an injective $D$-module because of the isomorphisms $hom_D(M, hom_A(D,K))\simeq hom_D(D{\otimes}_AM,K)\simeq hom_A(M,K)$. Indeed, with more details, if $0\rightarrow M' \stackrel{f}{\longrightarrow} M$ is a monomorphism of $D$-modules, applying $hom_A(\bullet,K)$ to this exact sequence, we get an epimorphism $hom_A(f,1_K)$ with the notation of [2] and thus an induced epimorphism $hom_D(f,1_{hom_A(D,K)})$. However, we emphasize once more that the left $D$-structure on $hom_A(D,K)$ used in [2] is coming from the right action of $D$ on $D=D_D$ through the formula $(\xi f)(P)=f(P\xi), \forall \xi\in T, \forall f\in hom_A(D,K)$ and therefore does not provide in general the structure of differential module defined by the formula $(\xi f)(P)=\xi(f(P))-f(\xi P)$ as in the theorem.  \\

\noindent
{\bf COROLLARY 3.16}: When $M\in mod(D)$ is any differential module, there is a sequence:  \\
\[   0 \rightarrow  {\wedge}^0T^*{\otimes}_AM \stackrel{\nabla}{\longrightarrow} {\wedge}^1T^*{\otimes}_AM \stackrel{\nabla}{\longrightarrow} ... \stackrel{\nabla}{\longrightarrow} {\wedge}^nT^*{\otimes}_AM \longrightarrow 0  \]
Accordingly, with $R=hom_A(M,A)={\wedge}^0T^*{\otimes}_A R$, there is a sequence:  \\
\[   0 \longrightarrow sol_A(M) \longrightarrow {\wedge}^0T^*{\otimes}_AR \stackrel{\nabla}{\longrightarrow} {\wedge}^1T^*{\otimes}_AR \stackrel{\nabla}{\longrightarrow} ... \stackrel{\nabla}{\longrightarrow} {\wedge}^nT^*{\otimes}_AR \longrightarrow 0     \]
where $sol_A(M)=hom_D(M,A)$. The corresponding deleted de Rham sequence $DR(R)$:\\
\[   0 \longrightarrow {\wedge}^0T^*{\otimes}_AR \stackrel{\nabla}{\longrightarrow} {\wedge}^1T^*{\otimes}_AR \stackrel{\nabla}{\longrightarrow} ... \stackrel{\nabla}{\longrightarrow} {\wedge}^nT^*{\otimes}_AR \longrightarrow 0  \]
only depends on the exterior derivative and the Spencer operator. \\

\noindent
{\it Proof}: With {\it any} differential module $M\in mod(D)$, we define an operator:  \\
\[    \nabla \hspace{2mm} :Ê\hspace{2mm} {\wedge}^sT^*{\otimes}_AM \rightarrow {\wedge}^{s+1}T^*{\otimes}_AM: {\wedge}^sT^*\otimes M_{q+r}\rightarrow {\wedge}^{s+1}T^*\otimes M_{q+r+1} \] 
extending the exterior derivative by the formula:  \\
\[   \nabla(\omega \otimes m)=d\omega\otimes m + dx^i\wedge \omega \otimes d_im, \hspace{3mm} \forall \omega\in {\wedge}^sT^*, \forall m\in M  \]
We obtain easily:\\
\[ {\nabla}^2(\omega\otimes m)= d^2\omega\otimes m + dx^i\wedge dx^j\wedge \omega \otimes d_{ij}m=0.  \]
The corresponding de Rham sequence $DR(M)$ may not be exact. Indeed, if $M=A$ with the canonical action of $P\in D$ on $a\in A$ given by $a\rightarrow P(a)$, then $DR(A)$ is just the Poincar\'{e} sequence and the first operator $\nabla=d$ is for sure not injective. On the contrary, when $M=D$, things are much more delicate and $DR(D)$ is exact unless at ${\wedge}^nT^*{\otimes}_AD$ where we have to add a surjective map ${\wedge}^nT^*{\otimes}_AD\rightarrow {\wedge}^nT^* :\alpha\otimes P\rightarrow \alpha. P:adx^1\wedge ... \wedge dx^n\rightarrow ad(P)(a)dx^1\wedge ... \wedge dx^n$ as we shall see in Theorem 3.23. \\
We also notice that $DR(M)=DR(A){\otimes}_AM$ as we are only concerned with the corresponding truncated sequences. Moreover, if $R$ is the inverse system of $M$, we may construct $DR(R)$ as before, replacing $m\in M$ by $f\in R$ and using the Spencer operator defined in the last theorem. Looking for the kernel of the first operator $f\rightarrow dx^i\otimes d_if$, we obtain $d_if=0\Leftrightarrow {\partial}_i(f(m))=f(d_im)\Leftrightarrow d_i\circ f=f\circ d_i, \forall i=1,...,n, \forall m\in M$ and thus $f\in hom_D(M,A)$ as claimed (See Example 3.20). \\
We now study the possibility to endow $DR(R)$ with a filtration as in the differential geometric setting of the first Spencer sequence, keeping in mind that, when $g_q$ is involutive ($2$-acyclic) and $g_{q+1}$ is a vector bundle over $X$, then $g_{q+r}$ is a vector bundle over $X$ for any $r\geq 0$($r\geq 1$) (See [14],[27] or [15] III.2.22,23 for more details). With $m=n=1, K=\mathbb{Q}(x)$, the counterexample $xy_x-y=0 \Rightarrow xy_{xx}=0$ is well known. If $G=gr(M)$ is the graded module of $M$ with $G={\oplus}^{\infty}_{q=0}G_q$, we have the short exact sequences $0\rightarrow M_{q} \rightarrow M_{q+1} \rightarrow G_{q+1} \rightarrow 0$ of modules over $A$ and it is tempting to compare them to the dual short exact sequences $0\rightarrow g_{q+1}\rightarrow R_{q+1}\rightarrow R_q \rightarrow 0$ that were used in the previous section. However, applying $hom_A(\bullet,A)$ to the first sequence does not in general provide a short exact sequence, unless the first sequence splits, that is if we replace vector bundles over $X$ used in Section 2 by finitely generated projective modules over $A$. One can also use a localization by introducing the field of fractions $K=Q(A)$ in order to deal only with finite dimensional vector spaces over $K$ or use the fact that $K=Q(A)$ is an injective module over $A$ and deal with $hom_A(\bullet,K)$ in order to obtain exact sequences. \\
In any case, $DR_{q+n+r}(R)$ starting with $R_{q+n+r}$ projects onto $DR_{q+n+r-1}(R)$ and the kernel of this projection is just (up to sign) the Spencer $\delta$-sequence:\\
\[  0\rightarrow g_{q+n+r}\stackrel{\delta}{\longrightarrow}T^*\otimes g_{q+n+r-1}\stackrel{\delta}{\longrightarrow}...\stackrel{\delta}{\longrightarrow}{\wedge}^nT^*\otimes g_{q+r}\rightarrow 0  \]
which is exact at the generic term $... \stackrel{\delta}{\longrightarrow} {\wedge}^sT^*\otimes g_{q+r} \stackrel{\delta}{\longrightarrow} ...$, $\forall r\geq 0$ whenever $g_q$ is involutive. However, as explained with details in ([14],[15]), a main formal problem is that the first Spencer sequence is {\it formally exact} (each operator generates the CC of the preceding one) but is {\it not strictly exact} (roughly the operators involved are not formally integrable and thus far from being involutive). Indeed, considering the system ${\partial}_if^k-f^k_i=0$ when $q=1$ and using crossed derivatives, we obtain the new first order equations ${\partial}_if^k_j-{\partial}_jf^k_i=0$. This is the reason for introducing the second Spencer sequence with Spencer bundles $C_r=C_{q,r}={\wedge}^rT^*\otimes R_q/\delta({\wedge}^{r-1}T^*\otimes g_{q+1})$ providing epimorphisms ${\wedge}^rT^*\otimes R_q\rightarrow C_{q,r}\rightarrow 0$ and $C_{q,r}\rightarrow {\wedge}^rT^*\otimes R_{q-1}\rightarrow 0$. Moreover, the second Spencer sequence being strictly exact with first order involutive operators, the sequence $\{C_{q+1,r}\}$ projects onto the sequence $\{C_{q,r}\}$ and we get the same projective limit as before if we replace vector bundles by projective modules over $A$.\\

It remains to compare the Spencer sequences thus obtained with the so-called Spencer sequence that can be found in the litterature on differential modules ([1],[7],[15],[26]), namely:  \\ 
\[   \begin{array}{rcccl}
0\longrightarrow D{\otimes}_AM {\otimes }_A{\wedge}^nT & \stackrel{d^*}{\longrightarrow}& D{\otimes}_AM {\otimes}_A{\wedge}^{n-1}T &\stackrel{d^*}{\longrightarrow}&   ...   \\
... &\stackrel{d^*}{\longrightarrow} &  D{\otimes}_AM{\otimes}_A{\wedge}^1T & \stackrel{d^*}{\longrightarrow}  & D{\otimes}_A M  \longrightarrow M\longrightarrow 0
\end{array}   \]
where we set, using the right structure over $D$ of the bimodule $D {\otimes}_A M$ and a hat for omission:\\
\[  \begin{array}{rcl}
d^*(P\otimes m \otimes{\xi}_1\wedge ... \wedge {\xi}_r) & = & {\sum}_i (-1)^{i-1}(P\otimes m){\xi}_i\otimes {\xi}_1\wedge ...\wedge{\hat{\xi}}_i\wedge ...\wedge {\xi}_r   \\
   & & +{\sum}_{i<j}(-1)^{i+j}P\otimes m \otimes [{\xi}_i,{\xi}_j]\wedge {\xi}_1\wedge ... \wedge {\hat{\xi}}_i\wedge ... \wedge {\hat{\xi}}_j\wedge ... \wedge {\xi}_r  
\end{array}   \]
Comparing to the standard definition of the exterior derivative, it is easy to check that $d^*\circ d^*=0$ while the last two maps are $P\otimes m \otimes \xi\rightarrow(P\otimes m)\xi=P\xi\otimes m-P\otimes \xi m$ and $P\otimes m \rightarrow Pm$. The left structure induced by $D={ }_DD$ commutes with $d^*$ and we have for example:  \\
\[  Qd^*(P\otimes \xi \otimes m)=QP\xi \otimes m -QP\otimes \xi m=d^*(QP\otimes \xi \otimes m)=d^*Q(P\otimes \xi\otimes m)\hspace{2mm} ,\forall \xi\in T, \forall P\in D  \]
or simply $Q((P\otimes m)\xi)=(Q(P\otimes m))\xi$ with $Q(P\otimes m)=QP\otimes m$. Meanwhile, the right structure induced by $D=D_D$ is used in order to act with $\xi$ on the right as in the last theorem. The corresponding deleted sequence is denoted by $SP(M)$ and we have $SP(M)=SP(A){\otimes}_AM$. It is important to notice that the construction of $SP(M)$ does not depend on any assumption made on the filtration of $M$. \\
We notice that $SP(M)$ is filtred by the complexes:  \\
\[  0 \rightarrow D_{q-n}{\otimes}_A M{\otimes}_A {\wedge}^nT \rightarrow ... \rightarrow D_{q-1}{\otimes}_A M{\otimes}_AT Örightarrow D_q{\otimes}_A M
\rightarrow M \rightarrow 0  \]
The associated graded complex is the tensor product by $M$ over $A$ of the {\it Koszul complex}: \\
\[  0 \rightarrow {\wedge}^nT{\otimes}_A S_{q-n}T \rightarrow  {\wedge}^{n-1}T{\otimes}_AS_{q-n+1}T \rightarrow ... \rightarrow T{\otimes}_AS_{q-1}T \rightarrow S_qT \rightarrow 0   \]
because we have the short exact sequences $0 \rightarrow D_{q-1} \rightarrow D_q \rightarrow S_qT \rightarrow 0 , \forall q>0$. This complex is thus exact because it is the dual of the trivial Spencer $\delta$-sequence and $M$ is locally free under the assumptions already made. Accordingly, $SP(M)$ is a resolution of $M$ but without any practical interest in view of the size of the modules invoved. \\
However, when $M$ is a filtred differential module, then $SP(M)$ is filtred by $SP_q(M)$, namely: \\ 
\[  \begin{array}{rcccl}
0\longrightarrow D{\otimes}_AM_{q-n} {\otimes }_A{\wedge}^nT & \stackrel{d^*}{\longrightarrow}& D{\otimes}_AM_{q-n+1} {\otimes}_A{\wedge}^{n-1}T &\stackrel{d^*}{\longrightarrow}&   ...   \\
... &\stackrel{d^*}{\longrightarrow} &  D{\otimes}_AM_{q-1}{\otimes}_A{\wedge}^1T & \stackrel{d^*}{\longrightarrow}  & D{\otimes}_A M_q  \longrightarrow M\longrightarrow 0
\end{array}  \]
and the construction is compatible because $D_1M_{q-r}\subset M_{q-r+1}$. The associated graded complex is the tensor product by $D$ of the dual of the Spencer $\delta$-sequence for the symbol and is thus exact, a result allowing to stabilize the cohomology of $SP_q(M)$ when $q$ is large enough and to get therefore a resolution of $M$ as its inductive limit is such a resolution. \\

When $E$ is a (finite dimensional) vector bundle over $X$/(finitely generated) projective module over $A$, we may apply the correspondence $J_{\infty}(E) \leftrightarrow D{\otimes}_AE^* : J_q(E)\leftrightarrow D_q{\otimes}_A E^*$ between jet bundles and {\it induced} left differential modules in order to be able to use the {\it double dual isomorphism} $E\simeq E^{**}$ in both cases. Hence, starting from a differential operator $E \stackrel{\cal{D}}{\longrightarrow}F$, we may obtain a finite presentation $D{\otimes}_AF \stackrel{\cal{D}}{\longrightarrow}D{\otimes}_AE \rightarrow M \rightarrow 0$ and conversely. We shall apply this procedure in order to study two particular cases before considering the general situation (Compare to [1], 1.5.1, p 29 where this motivating reference to differential geometry is lacking). Also, using the "local trivialty of projective modules", when $M$ is projective/locally free over $A$ ({\it care}), we recall that the functor $M{\otimes}_A \bullet$ is exact ([15], Proposition I.2.45, p 156; [25]). As we saw in section 2, this situation is realized, for example, when ${\pi}^{q+1}_q : R_{q+1} \rightarrow R_q$ is an epimorphism of vector bundles over $X$ and $g_q$ is involutive. In a word and contrary to Gr\"{o}bner algorithm, this is a modern version of the Janet algorithm for realizing the idea of Riqier to "cut"  the set of all derivatives of the unknowns into the "parametric" and "principal" subsets, under the condition that certain determinants should not vanish ([6], [23]).\\

When $M=A$ we may refer to the Poincar\'{e} sequence in section 2 and obtain the following strictly exact resolution of $A$ with induced left differential modules:  \\
\[     0\rightarrow D{\otimes}_A{\wedge}^nT\stackrel{d^*}{\longrightarrow}D{\otimes}_A{\wedge}^{n-1}T\stackrel{d^*}{\longrightarrow} ... \stackrel{d^*}{\longrightarrow}D{\otimes}_AT\stackrel{d^*}{\longrightarrow}D\rightarrow A\rightarrow 0 \]
because we have in this case $g_q=0, \forall q\geq 1 \Leftrightarrow G_q=0,\forall q\geq 1$ and thus $A=M_0=M_1=...=M$ as the gradient operator is a finite type trivially involutive first order operator. The operator $d^*$ is simply defined by the formula:  \\
\[  d^*(P \otimes d_{i_1}\wedge ... \wedge d_{i_s}) = {\sum}_r (-1)^{r-1}Pd_{i_r}\otimes d_{i_1}\wedge ...\wedge\hat{d_{i_r}}\wedge ...\wedge d_{i_s}   \]
This elementary example explains the need to introduce the bimodule $D$ on the left of the tensor products in order to obtain operators acting on the right as left $D$-linear maps and conversely by exchanging the words "left" and "right".\\

When $M=D$, the cokernel of the canonical inclusion $0 \rightarrow SP_{q-1}(D) \rightarrow SP_q(D)$ is the tensor product by $D$ over $A$ of the {\it Koszul complex} already provided. This complex is the dual of the Spencer $\delta$-sequence starting with $S_qT^*$ and is thus exact $\forall q>0$. As $SP_0(D)$ reduces to $0 \rightarrow D \stackrel{id_D}{\longrightarrow} D \rightarrow 0$ because $D{\otimes}_A A=D$ and $DD_q=D, \forall q\geq 0$, we can use an induction on $q$ showing that $SP(D)$ is a resolution of $D$. \\

In the general case, denoting by $H^s_q(M)$ the homology at ${\wedge}^sT{\otimes}_AM_{q-s}$ of $SP_q(M)$, we have to show that $H^s_q(M)=0, \forall q\gg 0$.  As we can always suppose that $M_q$ is involutive and defined by an involutive presentation of order $q$, then $DM_{q+r}=M, \forall r\geq 0$ and we just need to prove that $H^s_q(M)=0, \forall s\geq 0$ in order to prove that $SP_{q+r}(M)$ is a resolution of $M$, $\forall r\geq 0$ and thus that $SP(M)$ is a resolution of $M$. However, we have already exhibited {\it in this case} the strict short exact sequence $0 \rightarrow I \rightarrow D^l \rightarrow M \rightarrow 0$ providing the short exact sequences $0 \rightarrow I_{q+r} \rightarrow D^l_{q+r} \rightarrow M_{q+r} \rightarrow 0, \forall r\geq 0$. Introducing the short exact sequence $0 \rightarrow SP_{q}(I) \rightarrow SP_q(D^l) \rightarrow SP_q(M) \rightarrow 0$ of complexes, we deduce from the "{\it snake theorem}" in commutative algebra the following long exact Ò{\it connecting sequence}Ò ([2], p 29; [15], p 196-199; [25], p 171):\\
\[      ... \rightarrow H^s_q(D^l) \rightarrow H^s_q(M) \rightarrow H^{s-1}_q(I) \rightarrow ...    \]
However, we have $H^s_q(D^l)=0$ from the preceding particular case and we have $H^{s-1}_q(I)=0$ from the inductive assumption. It folows that $H^s_q(M)=0$ and the short exact sequence of complexes is in fact a resolution of the corresponding short exact sequence of differential modules, a result showing that $SP(M)$ is a resolution of $M$.  \\
Finally, acting component by component and using the adjoint isomoprphism for $mod(D)$ of the preceding theorem while changing the position of the terms in the tensor products, we get:  \\
\[  hom_D(SP(A){\otimes}_A M,A)\stackrel{\simeq}{\longrightarrow} hom_D(SP(A),hom_A(M,A)) \stackrel{\simeq}{\longrightarrow} DR(R)    \]   
because  $hom_D(D{\otimes}_A{\wedge}^rT,M)=hom_A({\wedge}^rT,M)={\wedge}^rT^*{\otimes}_A M$ as ${\wedge}^rT$ is free over $A$ with dimension $n!/(r!(n-r)!)$ and thus $ hom_D(SP(A),M)\simeq DR(M), \forall M\in mod(D)$. As a recapitulating comment, we may say that $SP(M)$ is made by left $D$-modules with left $D$-linear maps obtained by means of actions on the right while $DR(M{\otimes}_AD)$ is made by right $D$-modules with right $D$-linear maps obtained by means of actions on the left. Setting $n=r+s$, the conversion procedure that we shall recall in Theorem 3.23 allows to obtain:  \\
\[   hom_A({\wedge}^nT^*, {\wedge}^sT^*{\otimes}_AM_{q+s-n}{\otimes}_AD)={\wedge}^rT{\otimes}_AM_{q-r}{\otimes}_AD   \]
as a link between the two complexes.  \\

\hspace*{12cm}  Q.E.D.  \\

\noindent
{\bf REMARK 3.17}: When ${\cal{D}}=\Phi \circ j_q$ is an arbitrary but regular operator of order $q$, we may "{\it cut} " the Janet sequence at $F_0$ in two parts by introducing the systems $B_r= im({\rho}_r(\Phi))\subseteq J_r(F_0)$ with $B_0=F_0$ and $B_{r+1}\subseteq {\rho}_r(B_1)$ projecting onto $B_r, \forall r\geq 0$. When $\cal{D}$ is involutive, then $B_1\subseteq J_1(F_0)$ is also involutive with $B_{r+1}={\rho}_r(B_1), \forall r\geq 0$ and we have the "{\it truncated diagram} " linking the second Spencer sequence to a part of the Janet sequence:   \\
\[   \begin{array}{rcccccccl}
&0& &  0 &  &   &     &  0 &   \\
&\downarrow & & \downarrow  & & &  & \downarrow &  \\
&C_0  &\stackrel{D_1}{\longrightarrow} & C_1 & \stackrel{D_2}{\longrightarrow} & ...  & \stackrel{D_n}{\longrightarrow} & C_n & \rightarrow 0  \\
&\downarrow & & \downarrow & & &  & \downarrow & \\
0 \rightarrow E \stackrel{j_q}{\longrightarrow} & C_0(E)&\stackrel{D_1}{\longrightarrow} & C_1(E) & \stackrel{D_2}{\longrightarrow} & ... & \stackrel{D_n}{\longrightarrow} & C_n(E) & \rightarrow 0 \\
&\hspace{5mm} \downarrow {\Phi}_0& & \hspace{5mm}\downarrow {\Phi}_1 & &  & & \hspace{5mm} \downarrow {\Phi}_n & \\
&F_0  &\stackrel{{\cal{D}}_1}{\longrightarrow} & F_1 & \stackrel{{\cal{D}}_2}{\longrightarrow} &...  &  \stackrel{{\cal{D}}_n}{\longrightarrow} & F_n & \rightarrow 0 \\
&\downarrow  &  & \downarrow & & & & \downarrow &    \\
&0 &  & 0  & & && 0 &   
\end{array}   \]
where the epimorphisms ${\Phi}_1, ..., {\Phi}_n$ are successively induced by the epimorphism ${\Phi}_0=\Phi$, the canonical projection of $C_0(E)=J_q(E)$  onto $F_0=J_q(E)/R_q$ with $C_0=R_q$. It is known that the central sequence is locally exact. As we already pointed out that $g_q$ was a vector bundle, introducing the projection $R_{q-1}$ of $R_q$ into $J_{q-1}(E)$, we have $C_n\simeq {\wedge}^nT^*\otimes R_{q-1}$, $C_n(E)\simeq {\wedge}^nT^*\otimes J_{q-1}(E)$ and thus $F_n\simeq {\wedge}^nT^*\otimes (J_{q-1}(E)/R_{q-1})$. {\it It is not at all evident} that the dual of this diagram is nothing else but the resolution of the short exact sequence $0 \rightarrow I \rightarrow D^m \rightarrow M \rightarrow 0 $ considered in Proposition 3.6. Indeed, dualizing the diagram of Proposition 2.7, we obtain at once the following commutative and exact diagram: \\
\[   \begin{array}{cccccl}
       0   &   &   0   &   &   0   &   \\
    \uparrow &  & \uparrow &  & \uparrow  &   \\
 {\wedge}^{r-1}T{\otimes}_A G_{q+1}  & \stackrel{{\delta}^*}{ \leftarrow} & {\wedge}^rT{\otimes}_A G_q &
  \leftarrow & Z ( {\wedge}^rT{\otimes}_A G_q) &  \leftarrow 0  \\
      \parallel  &  &  \uparrow  &  & \uparrow  &  \\
 {\wedge}^{r-1}T{\otimes}_A G_{q+1} &\leftarrow  &{\wedge}^rT{\otimes}_A M_q & \leftarrow &C_r^* &
 \leftarrow 0  \\
      \uparrow  &  & \uparrow &  & \uparrow &   \\
        0  &  \leftarrow &  {\wedge}^rT{\otimes}_A M_{q-1} &  = & {\wedge}^r T{\otimes} M_{q-1} & \leftarrow 0  \\
         &   &  \uparrow &  &  \uparrow  &  \\
        &   &  0  &  &  0  &  
       \end{array}     \]
Applying the dual Spencer operator ${\wedge}^rT{\otimes}_AM_q \rightarrow {\wedge}^{r-1}T{\otimes}_AM_{q+1}$, we obtain the strictly exact 
{\it second Spencer sequence} $SSP_q(M)$:   \\
\[  0 \rightarrow D{\otimes}_AC_n^* \rightarrow D{\otimes}_AC_{n-1}^* \rightarrow ... \rightarrow D{\otimes}_AC_1^* \rightarrow D{\otimes }_AM_q \rightarrow M \rightarrow 0    \]
which is a resolution of $M$ {\it stabilizing the filtration at order} $q$ {\it only} by means of induced differential modules. Accordingly, the last two differential morphisms, induced by the morphisms $P\otimes \xi \otimes m \rightarrow P\xi \otimes m - P\otimes \xi m$ and $P \otimes m \rightarrow Pm$ of the sequence $... \rightarrow D{\otimes}_AT{\otimes}_AM_q\rightarrow D{\otimes}_AM_{q+1} \rightarrow M \rightarrow 0$, dualize the exact sequence $0 \rightarrow R_{q+r+1} \rightarrow J_{r+1}(C_0) \rightarrow J_r(C_1)$ as in ([15], p 367-369). \\

In the opinion of the author based on thirty years of explicit applications to mathematical physics (general relativity, gauge theory, theoretical mechanics, control theory), the differential geometric framework is quite more natural than the differential algebraic framework, the simplest example being the fact that the so-called {\it Cosserat equations} of elasticity theory, discovered by the brothers Eug\`{e}ne and Fran\c{c}ois Cosserat as early as in 1909, are nothing else but the formal adjoint $ad(D_1)$ of the first Spencer operator $D_1$ for the {\it Killing equations} in Riemannian geometry ([14], [15], [20], [21]). In particular, it must be noticed that the very specific properties of the Janet sequence, namely that it starts with an involutive operator of order $q\geq 1$ but the $n$ remaining involutive operators ${\cal{D}}_1, ..., {\cal{D}}_n$ are of order $1$ and in (reduced) Spencer form cannot be discovered from the differential module point of view. However, the importance of the torsion-free condition/test for differential modules is a novelty brought from the algebraic setting and known today to be a crucial tool for understanding control theory ([15]). Finally, the situation in the present days arrived to a kind of "{\it vicious circle} " because the study of differential modules is based on filtration and thus formal integrability while computer algebra is based on Gr\"{o}bner bases as a way to sudy the same questions but by means of highly non-intrinsic procedures as we saw.  \\

   We may compare the previous differential algebraic framework with its differential geometric counterpart. Indeed, using notations coherent with the ones of the previous section, if now ${\cal{D}}=\Phi \circ j_q:E \rightarrow F$ is an operator of order $q$ with $dim(E)=m,dim(F)=p$, we may consider the exact sequences $0 \rightarrow R_{q+r} \rightarrow J_{q+r}(E) \stackrel{{\rho}_r(\Phi)}{\longrightarrow} J_r(F)$ by introducing the $r$-prolongation of $\Phi$, induce the Spencer operator $D:R_{q+r+1}\rightarrow T^*\otimes R_{q+r}$ when $r\geq 0$ and pass to the projective limit $R=R_{\infty}$. In actual practice, when $r=1$ we have $a^{\tau\mu}_kf^k_{\mu}=g^{\tau} \Rightarrow a^{\tau\mu}_kf^k_{\mu+1_i}+({\partial}_i a^{\tau\mu}_k)f^k_{\mu}=g^{\tau}_i $ and thus $a^{\tau\mu}_k({\partial}_if^k_{\mu}-f^k_{\mu+1_i})={\partial}_ig^{\tau}-g^{\tau}_i$, a procedure that can be easily extended to any value of $r>0$. As a byproduct, the link existing with infinite jets can be understood by means of the following commutative and exact diagram:  \\
\[   \begin{array}{rccccc}
0 \rightarrow & R & \rightarrow & J(E) & \stackrel{\rho(\Phi)}{\longrightarrow} & J(F)  \\
                     &\hspace{3mm} \downarrow d  &     &  \downarrow d  &    &   \downarrow d     \\
0 \rightarrow & T^* \otimes R & \rightarrow & T^*\otimes J(E) & \stackrel{\rho(\Phi)}{\longrightarrow} & T^*\otimes J(F) 
\end{array}   \]
where $df=dx^i\otimes d_if$. Hence, using the Spencer operator on sections, we may characterize $R$ by the following equivalent properties (See [14], Proposition 10, p 83 for a nonlinear version that can be used in Example 2.26): \\

\noindent
1) $f\in R$ is killed by ${\rho}_r(\Phi)$ ({\it no differentiation of} $f$ {\it is involved} ), $\forall r\geq 0$.  \\
2) $f\in R \Rightarrow d_if\in R$ ({\it a differentiation of} $f$ {\it is involved} ), $\forall i=1,...,n$.  \\ 

As an equivalent differential geometric counterpart of the above result, we may also define the {\it r-prolongations} ${\rho}_r(R_q)=J_r(R_q)\cap J_{q+r}(E)$ of a given system $R_q\subset J_q(E)$ of order $q$ by applying successively the following formula involving the Spencer operator of the previous section:  \\
\[    {\rho}_1(R_q)=J_1(R_q)\cap J_{q+1}(E)=\{ f_{q+1}\in J_{q+1}(E) \mid f_q\in R_q, Df_{q+1}\in T^*\otimes R_q\}   \]
Now, if we have another system $R_{q+1}\subseteq {\rho}_1(R_q) \subset J_{q+1}(E)$ of order $q+1$ {\it and projecting onto} $R_q$, we have the commutative and exact diagram:\\
\[  \begin{array}{rcccl}
  & 0  &  & 0  &   \\
  & \downarrow &  & \downarrow &   \\
  0 \rightarrow & g_{q+1}  & \rightarrow  & {\rho}_1(g_q) &     \\
  & \downarrow & &  \downarrow &  \\
  0 \rightarrow & R_{q+1} & \rightarrow & {\rho}_1(R_q) &   \\
  &  \downarrow &  & \downarrow &   \\
  0 \rightarrow & R_q & = & R_q & \rightarrow 0  \\
  &  \downarrow &  & \downarrow  &  \\
  & 0 &  & 0 &    
  \end{array}     \]
 Chasing in this diagram, it follows that $R_{q+1}={\rho}_1(R_q)$ if and only if $g_{q+1}={\rho}_1(g_q)$ (Compare to Proposition 1.2.5 in [7]). Otherwise, we may start afresh with $R^{(1)}_q={\pi}^{q+1}_q(R_{q+1})$ (See Lemma III.2.46 in [15] for details).  \\

\noindent
{\bf REMARK 3.18}: As we shall see through the following examples, when a section $f\in R:y^k_{\mu} \rightarrow f^k_{\mu}\in A$ is given, it may not provide a formal power series solution. Accordingly, it may be useful to provide $f$ as a formal (in general infinite) summation $E\equiv f^k_{\mu}a^{\mu}_k=0$ called {\it modular equation} by Macaulay ([9], \S 59, p 67) and to set $d_iE\equiv ({\partial}_if^k_{\mu}-f^k_{\mu+1_i})a^{\mu}_k=0$. Equivalently, one can use ${\partial}_i$ on the coefficients of $E$ in $A$ and set $d_ia^{\mu}_k=0$ if ${\mu}_i=0$ or $d_ia^{\mu}_k=-a^{\mu-1_i}_k$ if ${\mu}_i>0$. When $A=K=k$ is a field of constants and $m=1$, we recover exactly the notation of Macaulay ({\it up to sign}) but the link with the Spencer operator has never been provided. With $n=3, m=1, q=2, K=\mathbb{Q}(x^1,x^2,x^3)$, the nice but quite tricky example $y_{33}-x^2 y_{11}=0, y_{22}=0$ provided by Janet (See [14] and [16] for more details) is such that $par=\{ y,y_1,y_2,y_3,y_{11}, y_{12},y_{13},y_{23},y_{111},y_{113}, y_{123},y_{1113}\}$ and can be generated by the unique modular equation $E\equiv a^{1113}+x^2a^{1333}+a^{12333}=0$ with $d_2E=0$, because $y_{12333}-y_{1113}=0, y_{1333}-x^2y_{1113}=0$ and all the jets of order $>5$ vanish (Exercise).\\

\noindent
{\bf EXAMPLE  3.19}: Coming back to Example 2.11 where $dim(g_2)=3, dim(g_3)=1, g_{4+r}=0, \forall r\geq 0$ and $g_3$ is $2$-acyclic, we have the following commutative and exact diagram where $dim(E)=1, dim(F_0)=3, dim(F_1)=18-15=3$:  \\
\[  \begin{array}{rcccccccl}
   & 0 & & 0  & & 0  &  & 0  &   \\
   & \downarrow  &  & \downarrow  & & \downarrow  & & \downarrow   &  \\
0 \rightarrow  & g_5 & \rightarrow  & S_5 T^* & \rightarrow & S_3T^*\otimes F_0 & \rightarrow & T^*\otimes F_1  & \rightarrow 0  \\
  &\hspace{2mm} \downarrow \delta &  &\hspace{2mm} \downarrow \delta & &\hspace{2mm} \downarrow \delta &  & \parallel &   \\
 0 \rightarrow &  T^*\otimes g_4  & \rightarrow& T^*\otimes S_4T^*& \rightarrow  & T^* \otimes S_2T^*\otimes F_0& \rightarrow &T^* \otimes F_1&   \rightarrow 0  \\
   &\hspace{2mm} \downarrow  \delta &   & \hspace{2mm} \downarrow \delta &  &\hspace{2mm} \downarrow \delta  &  & \downarrow  &  \\
 0 \rightarrow & {\wedge}^2T^*\otimes  g_3  &\rightarrow  & {\wedge}^2 T^*\otimes S_3 T^* & \rightarrow  &{\wedge}^2T^* \otimes T^* \otimes F_0  &\rightarrow& 0  \\
   & \hspace{2mm} \downarrow \delta &  &\hspace{2mm}  \downarrow \delta &  &\hspace{2mm}  \downarrow\delta  &   \\
 0 \rightarrow & {\wedge}^3 T^*\otimes g_2 & \rightarrow & {\wedge}^3T ^* \otimes S_2 T^*  & \rightarrow &{\wedge}^3T^*\otimes F_0& \rightarrow &0  \\
   & \downarrow &  & \downarrow &  & \downarrow     \\
   &  0  &    &  0 &    & 0 &   &  & 
  \end{array}   \]
and the long exact sequence where $dim(F_2)=28-45+18=1$:  \\
\[   0 \rightarrow S_6T^* \rightarrow S_4T^* \otimes F_0 \rightarrow S_2T^* \otimes F_1 \rightarrow F_2 \rightarrow 0   \]
providing the following free resolution with second order operators:  \\
\[   0 \rightarrow D  \rightarrow D^3 \rightarrow D^3 \rightarrow D   \rightarrow M  \rightarrow 0   \]
where the Euler-Poincar\'{e} characteristic is equal to $1-3+3-1=0$ as $M$ is defined by a finite type system. With a slight abuse of language while shifting the various filtrations, we may say that we have a strict resolution because all the operators involved, being homogeneous,  are formally integrable though not involutive. We finally notice that the first and second Spencer sequences coincide because we have $dim(R)=dim(R_3)=1+3+3+1=8$ as $g_4=0$.\\

\noindent
{\bf EXAMPLE 3.20}: With $n=1, m=1, q=2, A=\mathbb{Q}[x]\Rightarrow K=\mathbb{Q}(x)$ an thus $k=\mathbb{Q}$, let us consider the second order system $y_{xx}-xy=0$. We successively obtain by prolongation $y_{xxx}-xy_x-y=0, y_{xxxx}-2y_x-x^2y=0, y_{xxxxx}-x^2y_x-4xy=0, y_{xxxxxx}-6xy_x-(x^3+4)y=0$ and so on. We obtain the the corresponding board:Ê \\
\[  \begin{array}{r|c|c|c|c|c|c|c|l}
 order & y & y_x & y_{xx} & y_{xxx} & y_{xxxx} & y_{xxxxx} & y_{xxxxxx} &... \\
\hline
2 & -x & 0 & 1 & 0 & 0 & 0 & 0& ...   \\
3 & -1 & -x & 0 & 1 & 0 & 0& 0& ...   \\
4& -x^2 & -2  & 0 & 0 & 1 & 0& 0& ... \\
5& -4x & -x^2 & 0&0 & 0 & 1 & 0 & ... \\
6 & -(x^3+4) & -6x & 0 & 0 & 0 & 0 & 1 & ...
\end{array}  \]

Let us define the sections $f'$ and $f"$ by the following board where $d=d_x$: \\
\[  \begin{array}{r|c|c|c|c|c|c|c|l}
  section& y & y_x & y_{xx} & y_{xxx} & y_{xxxx} & y_{xxxxx} & y_{xxxxxx} &... \\
\hline
f' & 1 & 0 & x & 1 & x^2 & 4x & x^3+4& ...   \\
f" & 0 & 1 & 0 & x & 2 & x^2 & 6x& ...   \\
\hline
df' & 0 & -x  & 0 & -x^2 & -2x & - x^3 & -6x^2 & ... \\
df" & -1 & 0 & -x & -1 & -x^2 & -4x & -x^3-4 & ... \\
\end{array}  \]
in order to obtain $df'=-xf", df"=-f'$. Though this is not evident at first sight, {\it the two boards are orthogonal over} $K$ in the sense that each row of one board contracts to zero with each row of the other though only the rows of the first board do contain a finite number of nonzero elements. {\it It is absolutely essential to notice that the sections} $f'$ {\it and} $f"$ {\it have nothing to do with solutions}  because $df'\neq 0, df"\neq 0$ on one side and also because $d^2f'-xf'=-f"=\frac{1}{x}df'\neq 0$ even though $d^2f"-xf"=0$ on the other side. As a byproduct, $f'$ or $f"$ can be chosen separately as unique generating section of the inverse system over $K$ ({\it care}) and we may write for example $f' \rightarrow E'\equiv a^0+xa^{xx}+a^{xxx}+x^2 a^{xxxx}+ ... =0$ while $f" \rightarrow E"\equiv a^x+xa^{xxx}+2a^{xxxx}+ ... =0$. \\
Finally, setting $f=af'+bf"$, we have $ df=(\partial a)f'+(\partial b-xa)f"=0\Leftrightarrow {\partial}^2a-xa=0, b=\partial a$. If $a=P/Q$ with $P,Q\in \mathbb{Q}[x]$ and $Q\neq 0$, we obtain easily :Ê\\
\[  Q^2{\partial}^2P-2Q\partial P\partial Q-PQ {\partial}^2Q +  2P(\partial Q)^2-xPQ^2=0  \]
If $deg(P)=p, deg(Q)=q$, the four terms on the left have the same degree $p+2q-2$ while the last term has degree $p+2q+1$ and thus $Q\neq 0 \Rightarrow P=0 \Rightarrow a=0 \Rightarrow b=0$.  \\
With $D=A[d], \Phi\equiv y_{xx}-xy=z \Rightarrow d_x\Phi \equiv y_{xxx}-xy_x-y=z_x, ... $ and ${\cal{D}}=\Phi \circ j_2$, we have the finite presentation $0 \rightarrow D \stackrel{\cal{D}}{\longrightarrow} D \rightarrow M \rightarrow 0$. Using sections with $f_{xx}-xf=g, f_{xxx}-xf_x-f=g_x, ...$ and $g_r=j_r(x)$, we may choose $f=-1, f_x=0, f_{xx}=0, ...$. However, in a similar situation with $\Phi \equiv xy_x-y \Rightarrow d_x \Phi \equiv xy_{xx}, ...$ and $g_r=j_r(x)$, we may choose $f=0, f_x=1$ but obtain $xf_{xx}=1$ which {\it cannot} be solved over $A$, an example showing why the divisibility of $K$ over $A$ is needed in order to get the short exact sequence $0 \rightarrow hom_A(M,K) \rightarrow hom_A(D,K) \rightarrow hom_A(D,K) \rightarrow 0$. \\

\noindent
{\bf EXAMPLE 3.21}: With $n=1, m=2, q= 2, k=\mathbb{Q}$, let us consider the system $y^1_{xx}=0,y^2_x=0$. Setting  
$z^1=y^1, z^2=y^1_x, z^3=y^2$, we obtain the involutive system:
\[  \left\{   \begin{array}{lcc}
z^1_x-z^2 & = & 0 \\
z^2_x& =  & 0 \\
z^3_x  & = & 0 
\end{array}
\right. \fbox{$\begin{array}{c}
 x \\
 x  \\
 x
 \end{array}$}  \]
Setting $d=d_x$, then $\mathfrak{m}=(d)$ kills $y^1_x$ {\it and} $y^2$, that is $z^2$ {\it and} $z^3$. It follows that the CK data for $z=g(x)$ are 
$\{g^1(x)=g^1(0), g^2(x)=g^2(0), g^3(x)=g^3(0)\}$ and we have the finite basis: \\
\[  \begin{array}{l|ccc|ccc|ccc|c}
    & z^1 & z^2 & z^3 & z^1_x & z^2_x & z^3_x & z^1_{xx} &z^2_{xx} & z^3_{xx} & ... \\
\hline
g  & 1 & 0 & 0 & 0 & 0 & 0 & 0 & 0 & 0 & ... \\
g' & 0 & 1 & 0 & 1 & 0 & 0 & 0 & 0 & 0 & ... \\
g" & 0 & 0 & 1 & 0 &0 & 0 & 0 & 0 & 0 & ... 
\end{array}  \] 
As $dg=0,dg'=-g, dg"=0$, a basis with only two generators may be $\{g',g"\}$. However, if we work with the differential field $K=k(x)$ instead of $k$, we may introduce $y=y^2-xy^1$ which must satisfy $y_{xxx}=0$ over $k$ and there is a unique generator according to ([19], Ex. 3.8).   \\ 
Similarly, with the system $y^1_{xx}-y^1=0, y^2_x=0$, we may consider the involutive system:\\
\[  \left\{   \begin{array}{lcc}
z^1_x-z^2 & = & 0 \\
z^2_x-z^1& =  & 0 \\
z^3_x  & = 0 
\end{array}
\right. \fbox{$\begin{array}{c}
 x \\
 x  \\
 x
 \end{array}$}  \]
and consider now:  \\
\[  \begin{array}{l|ccc|ccc|ccc|c}
    & z^1 & z^2 & z^3 & z^1_x & z^2_x & z^3_x & z^1_{xx} & z^2_{xx} & z^3_{xx} & ... \\
\hline
g  & 1 & 0 & 0 & 0 & 1 & 0 & 1 & 0 & 0 & ... \\
g' & 0 & 1 & 0 & 1 & 0 & 0 & 0 & 1 & 0 & ... \\
g" & 0 & 0 & 1 & 0 &0 & 0 & 0 & 0 & 0 & ... 
\end{array}  \] 
We obtain $dg=-g', dg'=-g, dg"=0$ and thus:  \\
\[   h=g-g", dh=-g', d^2h=g  \Leftrightarrow g=d^2h, g'=-dh, g"=d^2h-h  \]
that is a unique generator $h$ as in ([9],$\S 72$, p 81). We may also introduce $y=y^2-y^1$ which must satisfy $y_{xxx}-y_x=0$ and there is a unique generator as in ([19], Ex. 3.8). \\

\noindent
{\bf EXAMPLE 3.22}: With $n=3, m=3, p=2, q=2, k=\mathbb{Q}$, let us consider the differential module $M$ defined by the second order system: \\ 
\[  \left\{  \begin{array}{lcl}
y^3_{33}-y^1_1 & = & 0 \\
y^3_{23}-y^2_1-y^3 & = 0 
\end{array}
\right. \fbox{$\begin{array}{lll}
1 & 2 & 3 \\
1 & 2 & \bullet
\end{array} $}  \]
Using the only non-multiplicative variable involved in the board, the system is not formally integrable and we may consider anew the second order 
system: \\
\[  \left\{  \begin{array}{lcl}
y^3_{33}-y^1_1 & = & 0 \\
y^3_{23}-y^2_1-y^3 & = & 0 \\
y^2_{13}-y^1_{12}+y^3_3 & = & 0
\end{array}
\right. \fbox{$\begin{array}{lll}
1 & 2 & 3 \\
1 & 2 & \bullet \\
1 & \bullet & \bullet
\end{array} $}  \]
The coordinate system is not $\delta$-regular and we may change the coordinates with $x^1\rightarrow x^1, x^2 \rightarrow x^2 , x^3 \rightarrow x^3+x^1$ in order to get the involutive system with ${\alpha}^3_2=1$:  \\
\[  \left\{  \begin{array}{lcl}
y^3_{33}-y^1_3 -y^1_1 & = & 0 \\
y^2_{33}+ y^2_{13}-y^1_{23}-y^1_{12}+y^3_3 & = & 0 \\
y^3_{23}-y^2_3-y^2_1-y^3 & = & 0
\end{array}
\right. \fbox{$\begin{array}{lll}
1 & 2 & 3 \\
1 & 2 & 3 \\
1 & 2 & \bullet
\end{array} $}  \]
We obtain therefore the strict free resolution $0 \longrightarrow D \longrightarrow D^3 \longrightarrow D^3 \longrightarrow M \longrightarrow 0 $. \\
Localizing the initial system or using computer algebra, we obtain:  \\
\[  ({\chi}_3)^2 y^3-{\chi}_1y^1=0, ({\chi}_2{\chi}_3-1)y^3-{\chi}_1y^2=0 \Rightarrow {\chi}_1 (({\chi}_3)^2y^2-({\chi}_2{\chi}_3-1)y^1=0 \]
Accordingly, $z=y^2_{33}-y^1_{23}+y^1$ generates $t(M)$ as we get $d_1 z=z_1=0$. The torsion-free module $M'=M/t(M)$ is defined, in the initial coordinate system, by the second order system: \\
\[  \left\{  \begin{array}{lcl}
y^3_{33}-y^1_1 & = & 0 \\
y^2_{33}-y^1_{23}+y^1 & = & 0 \\
y^3_{23}-y^2_1-y^3 & = & 0 \\
y^2_{13}-y^1_{12}+y^3_3 & = & 0
\end{array}
\right. \fbox{$\begin{array}{lll}
1 & 2 & 3 \\
1 & 2 & 3 \\
1 & 2 & \bullet \\
1 & \bullet & \bullet
\end{array} $}  \]
The coordinate system is not $\delta$-regular and we may change the coordinates with $x^1 \rightarrow x^1, x^2 \rightarrow x^2+x^1, x^3 \rightarrow x^3 $ in order to obtain the second order involutive system:  \\
\[  \left\{  \begin{array}{lcl}
y^3_{33}-y^1_2-y^1_1  & = & 0 \\
y^2_{33}- y^1_{23}+y^1 & =  &0 \\
y^3_{23}-y^2_2 - y^2_1- y^3 & = & 0 \\
y^2_{23}+ y^2_{13}-y^1_{22}-y^1_{12}+ y^3_3 & = & 0  
\end{array}
\right. \fbox{$\begin{array}{lll}
1 & 2 & 3 \\
1 & 2 & 3 \\
1 & 2 & \bullet \\
1 & 2 & \bullet
\end{array} $}  \]
and the strict free resolution $0 \longrightarrow D^2 \longrightarrow D^4 \longrightarrow D^3 \longrightarrow M' \longrightarrow 0$.  \\
Using double duality, one can exhibit the injective parametrization ${\xi}_{33}=y^1, {\xi}_{23}-\xi=y^2, {\xi}_1=y^3$ in the original coordinate system (exercise) and thus $M'\simeq D$ is free. Also, using the solution of the {\it Bose conjecture} ([16],p 216-219), one can define $M'$ by the only two PD equations (exercise):  \\
\[  \left\{  \begin{array}{lcl}
y^2_{33}-y^1_{23}+ y^1  & = & 0 \\
y^3 + y^2_{123}-y^1_{122}+ y^2_1 & = & 0 
\end{array}
\right.   \]
As the second PD equation provides $y^3$, there is no CC and we get the free resolution $0 \longrightarrow D^2 \longrightarrow D^3 \longrightarrow M' \longrightarrow 0 $ which is however not strict (See [19] for more details on strict morphisms and resolutions). Finally, integrating by part ${\lambda}^2(y^3_{33}-y^1_1)+{\lambda}^1(y^3_{23}-y^2_1-y^3)$ and permuting the coordinates with $(123)\rightarrow (312)$, we get the adjoint system ${\lambda}^2_3=0, {\lambda}^1_3=0, {\lambda}^2_{22}+{\lambda}^1_{12}-{\lambda}^1=0$ leading to ${\tilde{\alpha}}^3_2=0$ and we check that $m-{\alpha}^3_2=3-1=2=p-{\tilde{\alpha}}^3_2$, {\it a general result not evident at all}.   \\

As $D={ }_DD_D$ is a bimodule, then $M^*=hom_D(M,D)$ is a right $D$-module according to Lemma 3.1 and we may thus define a right module $N_r$ by the ker/coker long exact sequence $0\longleftarrow N_r \longleftarrow F_1^*\stackrel{{\cal{D}}^*}{ \longleftarrow} F^*_0 \longleftarrow M^* \longleftarrow 0$.\\

\noindent
{\bf THEOREM 3.23}: We have the {\it side changing} procedures $M={ }_DM \rightarrow M_r={\wedge}^nT^*{\otimes}_AM$ and $N=N_D \rightarrow N_l=hom_A({\wedge}^nT^*,N_r)$ with $(M_r)_l=M$ and $(N_l)_r=N$. \\

\noindent
{\it Proof}: According to the above Theorem, we just need to prove that ${\wedge}^nT^*$ has a natural right module structure over $D$. For this, if $\alpha=adx^1\wedge ...\wedge dx^n\in {\wedge}^nT^*$ is a volume form with coefficient $a\in A$, we may set $\alpha.P=ad(P)(a)dx^1\wedge...\wedge dx^n$ when $P\in D$. As $D$ is generated by $A$ and $T$, we just need to check that the above formula has an intrinsic meaning for any $\xi={\xi}^id_i\in T$. In that case, we check at once:
\[  \alpha.\xi=-{\partial}_i(a{\xi}^i)dx^1\wedge...\wedge dx^n=-\cal{L}(\xi)\alpha \]
by introducing the Lie derivative of $\alpha$ with respect to $\xi$, along the intrinsic formula ${\cal{L}}(\xi)=i(\xi)d+di(\xi)$ where $i( )$ is the interior multiplication and $d$ is the exterior derivative of exterior forms. According to well known properties of the Lie derivative, we get :
\[\alpha.(a\xi)=(\alpha.\xi).a-\alpha.\xi(a), \hspace{5mm} \alpha.(\xi\eta-\eta\xi)=-[\cal{L}(\xi),\cal{L}(\eta)]\alpha=-\cal{L}([\xi,\eta])\alpha=\alpha.[\xi,\eta].  \]
Using the anti-isomorphism $ad:D \rightarrow D : P \rightarrow ad(P)$, we may also introduce the {\it adjoint functor} $ad : mod(D) \rightarrow mod (D^{op}): M \rightarrow ad(M)$ with $for(M)=for(ad(M))$ and $m.P=ad(P)m, \forall m\in M, \forall P\in D$. We obtain:  \\
\[    m.(PQ)=ad(PQ)m=(ad(Q)ad(P))m=ad(Q)(ad(P)m)=(m.P).Q, \forall P,Q\in D  \]
We have an $A$-linear isomorphism $ad(M)\simeq M_r:m \rightarrow m\otimes \alpha$ in $mod(D^{op})$. Indeed, with $\alpha=dx^1\wedge ... \wedge dx^n$ and any $d$ among the $d_i$ in place of $\xi$, we get $m.d=ad(d)m=-dm$ in $ad(M)$ while $(m\otimes \alpha)d=-dm\otimes \alpha, \forall m\in M$ in $M_r$ because $div(d)=0$ and thus ${\cal{L}}(d)\alpha =0$. Accordingly, the previous isomorphism is right $D$-linear.  \\
In order to study the case of $D={ }_DD$, considered as a left $D$-module over $D$, we shall compare $ad(D), D_r$ and $D_D$. According to the last isomorphism obtained, we just need to study the isomorphim $ad(D)\simeq D_D: P \rightarrow ad(P)$. Indeed, we get $P \rightarrow P.Q=ad(Q)P\neq PQ$ and obtain therefore $P.Q \rightarrow ad(P.Q)=ad(ad(Q)P)=ad(P)ad(ad(Q))=ad(P)Q$, a result showing that this isomorphism is also right $D$-linear. \\
\hspace*{12cm}  Q.E.D.  \\
 
\noindent
{\bf REMARK 3.24}: The above results provide a new light on duality in physics. Indeed, as the Poincar\'{e} sequence is self-adjoint (up to sign) 
{\it as a whole} and the linear Spencer sequence for a system of finite type is locally isomorphic to copies of that sequence, it follows {\it in this case} from Proposition 3.4 that $ad(D_{r+1})$ parametrizes $ad(D_r)$ in the dual of the Spencer sequence while $ad({\cal{D}}_{r+1})$ parametrizes $ad({\cal{D}}_r)$ in the dual of the Janet sequence, {\it a result highly not evident at first sight because} ${\cal{D}}_r$ {\it and} $D_{r+1}$ {\it are totally different operators}. The reader may look at [20] and [21] for recent applications to mathematical physics (gauge theory and general relativity). \\
 
We shall now study with more details the module $M$ versus the system $R$ when $D=K[d]$. First of all, as $K$ is a field, we obtain in particular the Hilbert polynomial $dim_K(M_{q+r})=dim_K(R_{q+r})= \frac{\alpha}{d!}r^d+ ...$ where the intrinsic integer $\alpha$ is called the {\it multiplicity} of $M$ and is the smallest non-zero character. We use to set $d(M)=d \Rightarrow cd_D(M)=cd(M)=n-d=r, rk_D(M)=rk(M)=\alpha$ if $cd(M)=0$ and $0$ otherwise. Now, If $M$ is a module over $D$ and $m\in M$, then the cyclic differential submodule $Dm\subset M$ is defined by a system of OD or PD equations for one unknown and we may look for its codimension $cd(Dm)$. A similar comment can be done for any differential submodule $M'\subset M$. Sometimes, a single element $m\in M$, called {\it differentially primitive element}, may generate $M$ if $Dm=M$ as in Example 3.15. Using the results of ([7]), we get:    \\
 
\noindent
{\bf PROPOSITION 3.25}: $t_r(M)=\{m\in M \: {\mid} \: cd(Dm)> r \}$ is the greatest differential submodule of $M$ having codimension $> r$ and does not depend on the presentation or filtration of $M$.  \\

\noindent
{\bf PROPOSITION 3.26}: $cd(M)=r \Longleftrightarrow {\alpha}^{n-r}_q\neq 0, {\alpha}^{n-r+1}_q= ... ={\alpha}^n_q=0 \Longleftrightarrow t_r(M)\neq M, t_{r-1}(M)= ... =t_0(M)=t(M)=M$ and this intrinsic result can be most easily checked by using the standard or reduced Spencer form of the system defining $M$.  \\

We are now in a good position for defining and studying {\it pure differential modules} along lines similar to the ones followed by Macaulay in ([9]) for studying 
{\it unmixed polynomial ideals}. \\

\noindent
{\bf DEFINITION 3.27}: $M$ is $r$-{\it pure} $\Longleftrightarrow t_r(M)=0, t_{r-1}(M)=M \Longleftrightarrow cd(Dm)=r, \forall m\in M$. More generally, $M$ is {\it pure} if it is $r$-pure for a certain $0\leq r\leq n$ and $M$ is {\it pure} if it is $r$-pure for a certain $0\leq r \leq n$. In particular, $M$ is $0$-pure if $t(M)=0$ and, if $cd(M)=r$ but $M$ is not $r$-pure, we may call $M/t_r(M)$ the {\it pure part} of $M$. It follows that $t_{r-1}(M)/t_r(M)$ is equal to zero or is 
$r$-pure (See the picture in [15], p 545). When $M=t_{n-1}(M)$ is $n$-pure, its defining system is a finite dimensional vector space over $K$ with a symbol of finite type (see Example 2.11 of 2-acyclicity). Finally, when $t_{r-1}(M)=t_r(M)$, we shall say that there is a {\it gap} in the purity filtration:   \\
\[   0=t_n(M) \subseteq t_{n-1}(M) \subseteq ... \subseteq t_1(M) \subseteq t_0(M)=t(M) \subseteq M     \]

\noindent
{\bf EXAMPLE 3.28}: If $K=\mathbb{Q}$ and $M$ is defined by the involutive system $y_{33}=0, y_{23}=0, y_{13}=0$, then $z=y_3$ satifies $d_3z=0, d_2z=0, d_1z=0$ and $cd(Dz)=3$ while $z'=y_2$ only satisfies $d_3z'=0$ and $cd(Dz')=1$. We have the purity filtration 
$0 =t_3(M) \subset t_2(M) =t_1(M) \subset t_0(M)=t(M)=M$ with one gap and two strict inclusions.  \\

When $k$ is a field of constants, then $D=k[d]$ is isomorphic to the polynomial ring $A=k[\chi]$ and we shall generalize the criterion of Macaulay in ([9], $\S 41$, p. 43) as follows when $cd(M)=r$ by setting $\chi=({\chi}',{\chi}")$ with ${\chi}'=({\chi}_1,...,{\chi}_{n-r})$ and ${\chi}"=({\chi}_{n-r+1},...,{\chi}_n)$ while introducing similarly $d=(d',d")$ with $d'=(d_1,...,d_{n-r})$ and $d"=(d_{n-r+1},...,d_n)$ or $\mu=({\mu}',{\mu}")$ with ${\mu}'=({\mu}_1,...,{\mu}_{n-r})$ and ${\mu}"=({\mu}_{n-r+1},...,{\mu}_n)$. For such a purpose we recall a few results about the localization used in the {\it primary decomposition} of a module $M$ over a commutative integral domain $A$ which are not so well known ([15],[18],[19]). We denote as usual by $spec(A)$ the set of {\it prime ideals} in $A$, by $max(A)$ the subset of {\it maximal ideals} in $A$ and by $ass(M)=\{\mathfrak{p}\in spec(A) {\mid} \exists 0\neq m\in M, \mathfrak{p}=ann_A(m)\}$ the (finite) set $\{ {\mathfrak{p}}_1, ... , {\mathfrak{p}}_t\}$ of {\it associated prime ideals}, while we denote by $\{ {\mathfrak{p}}_1, ...{\mathfrak{p}}_s\}$ the subset of {\it minimum associated prime ideals}. It is well known that $M\neq 0 \Longrightarrow ass(M)\neq \emptyset$. We recall that an ideal $\mathfrak{q}\subset A$ is 
$\mathfrak{p}$-{\it primary} if $ab\in \mathfrak{q}, b\notin \mathfrak{q} \Longrightarrow a\in rad(\mathfrak{q})=\mathfrak{p}\in spec(A)$. We say that a module $Q$ is $\mathfrak{p}$-{\it primary} if $am=0, 0\neq m\in Q \Longrightarrow a\in \mathfrak{p}=rad(\mathfrak{q})\in spec(A)$ when $\mathfrak{q}=ann_A(Q)$ or, equivalently, $ass(Q)=\{\mathfrak{p}\}$. Similarly, we say that a module $P$ is $\mathfrak{p}$-{\it prime} if $am=0, 0\neq m\in P \Longrightarrow a\in \mathfrak{p}\in spec(A)$ when $\mathfrak{p}=ann_A(P)$. It follows that any $\mathfrak{p}$-prime or $\mathfrak{p}$-primary module is $r$-pure with $n-r=trd(A/\mathfrak{p})$. Accordingly, a module $M$ is $r$-pure if and only if $\mathfrak{a}=ann_A(M)$ admits a primary decomposition $\mathfrak{a}={\mathfrak{q}}_1\cap ... \cap {\mathfrak{q}}_s$ and $rad(\mathfrak{a})={\mathfrak{p}}_1\cap ... \cap {\mathfrak{p}}_s$ with $cd(A/{\mathfrak{p}}_i)=cd(M)=r, \forall i=1, ... ,s$. {\it In that case only}, the monomorphism $0 \longrightarrow M \longrightarrow {\oplus}_{\mathfrak{p}\in ass(M)}M_{\mathfrak{p}}$ induces a monomorphism $0 \longrightarrow M \longrightarrow Q_1\oplus ... \oplus Q_s$ called {\it primary embedding} where the primary modules $Q_i$ are the images of the localization morphisms $M \longrightarrow M_{{\mathfrak{p}}_i}=S_i^{-1}M$ with $S_i=A-{\mathfrak{p}}_i$ inducing epimorphisms $M \longrightarrow Q_i \longrightarrow 0$ for $i=1, ... ,s$. Macaulay was only considering the case $M=A/\mathfrak{a}$ with a primary decomposition $\mathfrak{a}={\mathfrak{q}}_1\cap ... \cap {\mathfrak{q}}_s$.  \\

\noindent
{\bf THEOREM 3.29}: $M$ is $r$-pure {\it if and only if } $S({\chi}')z\in I$ {\it requires} $z\in I$ that is a residue $\bar{z}=0$ in $M$.  \\

\noindent
{\it Proof}: Let $I\subset F$ be the module of equations of $M=F/I$ and consider a primary decomposition $I={\cap}I_i$ of $I$ in $F$. We may pass to the residue by introducing short exact sequences $0\rightarrow I_i\rightarrow F \rightarrow Q_i \rightarrow 0$ providing induced epimorphisms $M\rightarrow Q_i\rightarrow 0$ both with a monomorphism $0\rightarrow M\rightarrow {\oplus}_iQ_i$ called {\it primary embedding} ([15], p 113). The $Q_i$ are primary modules and we may introduce the primary ideals ${\mathfrak{q}}_i=ann(Q_i)$ with ${\mathfrak{p}}_i=rad({\mathfrak{q}}_i)\in spec(A)$ in order to obtain the primary decomposition $\mathfrak{a}={\cap}{\mathfrak{q}}_i$ as in ([15], p 112). \\
If $M$ is not $r$-pure, then $\mathfrak{a}$ is mixed and we may suppose that $cd(A/{\mathfrak{p}}_t)>r$. Then we may choose $S\in {\mathfrak{q}}_t$ and $z\notin I_t$ while $z\in I_i, S\in A \Rightarrow Sz\in I_i, \forall i< t$ and thus $Sz\in {\cap}_{i< t}I_i$, that is to say $Sz\in {\cap}I_i=I$ though $z\notin I$ because $z\notin I_v$. Conversely, if $M$ is $r$-pure, then $\mathfrak{a}$ is unmixed with $s=t$ {\it necessarily}, each associated prime has codimension $r$, no one can thus contain $S$ and $Sz\in I\Rightarrow z\in I$ from the next argument. \\
Let us prove that, if $\mathfrak{b}\subset A$ is an ideal such that $I\subset I:\mathfrak{b}=J\neq I$ in $F$, then $\mathfrak{b}\subset {\mathfrak{p}}_i$ for a certain $i$. Indeed, if $z\in J\subset F,z\notin I$, then $bz\in I,\forall b\in \mathfrak{b}$ and thus $b\bar{z}=0$ with $\bar{z}\neq 0$ in $M$. Accordingly, $b$ is a zero-divisor in $M$ and thus $b\in \cup {\mathfrak{p}}_i$ the set of all zero-divisors of $M$ ([15], p 101). Equivalently, $z\notin I_i$ for a certain $i$ because otherwise we should have $z\in \cap I_i=I$.But $I_i$ is primary in $F$ with $bz\in I\subset I_i, z\notin I_i$ and thus $b\in {\mathfrak{p}}_i\Rightarrow \mathfrak{b}\subset {\mathfrak{p}}_i$ in any case. Hence, if $I:\mathfrak{m}\neq I$ for $\mathfrak{m}\in max(A)$, then $\mathfrak{m}={\mathfrak{p}}_i$ for a certain $i$ and thus $\mathfrak{m}\in ass(M) $.  \\
\hspace*{12cm}    Q.E.D.   \\

Using standard localization techniques, we get at once (See [19] for more details):  \\

\noindent
{\bf COROLLARY 3.30}: There is an exact sequence $0 \rightarrow t_r(M) \rightarrow M \rightarrow k({\chi}')\otimes M$ and the differential module $M$ is thus $r$-pure if and only if  $cd(M)=r$ and there is a monomorphism $0 \longrightarrow M \longrightarrow k({\chi}') \otimes M  $. In that case, $k({\chi}')$ is playing the part of $k$ in the localization of $M$ which is finite type over $k({\chi}')[d"]$. If $M$ is defined over $D=K[d]$, then $M$ is $r$-pure if and only if 
$cd(M)=r$ (the classes $n-r+1,...,n$ are full in any presentation) and the differential module defined by the equations of class $1$+...+class $(n-r)$ of a Spencer form is torsion-free.  \\

When $1\in S$ is a multiplicative subset of an integral domain $A$ and $M$ is a module over $A$, we have $S^{-1}M=S^{-1}A{\otimes}_AM$. When $M$ is finitely presented over $A$, we obtain for any other module $N$ over $A$ the so-called {\it localizing isomorphism} as in ([15], p 153) or( [24], Th 3.84,p. 107)), namely:  \\
\[  S^{-1}hom _A(M,N) \simeq hom_A(M,S^{-1}N) \simeq hom_{S^{-1}A}(S^{-1}M,S^{-1}N)   \]
by introducing $S^{-1}f:S^{-1}M \rightarrow S^{-1}N: s^{-1}x \rightarrow s^{-1}y$ over $S^{-1}A$ whenever $f:M \rightarrow N : x \rightarrow y$ with $y=f(x)$ is given over $A$. \\
As $M$ is finitely presented, we obtain with $S=k[{\chi}']\subset k[\chi]$ the isomorphism:Ê \\\[hom_{k({\chi}')}(k({\chi}')\otimes M,k({\chi}'))\simeq k({\chi}')\otimes hom_k(M,k)=k({\chi}')\otimes R\]
providing the localized system over $k({\chi}')[d"]$, on the condition to consider $k' =k({\chi}')$ as a new field of constants for $d"$ along lines first proposed by Macaulay ([9],$\S 43,77$) in order to extend maps $M \rightarrow k$ to maps $k({\chi}')\otimes M  \rightarrow k({\chi}')$, that is to 
maps $k' \otimes MÊ\rightarrow k' $. Moreover, the new differential module $M"=M/t_r(M)$ is $r$-pure {\it in any case} and provides a subsystem $R"\subseteq R$.  \\

\noindent
{\bf COROLLARY 3.31}: The differential module $M$ is $r$-pure if and only if  $cd(M)=r$ and there is an isomorphism $R"\simeq R$.  \\

As the first order CC of an involutive system are {\it automatically} in reduced Spencer form, we may just use the column $n-r$ of the Janet board in order to discover that the above {\it relative localization} " kills "  the equations of class $1$ up to class $n-r-1$ as we can divide by ${\chi}_{n-r}$ and thus only depends on the equations of strict class $n-r$ providing the nonzero intrinsic character  $\alpha$ when $r<n$ (Compare to [9], $\S 43$ and see [19], Proposition 5.7 for more details). Accordingly, as the value of this intrinsic character $\alpha$ may be computed for a Spencer form defining $M$ and we look for the dimension of $k'\otimes M$ as a vector space over $k'$ which is equal to $\alpha$ when $k'\neq k$, we obtain the following rather striking result which could not even be imagined without involution:   \\

\noindent
{\bf COROLLARY 3.32}: The localized system is finite type with a (finite) dimension equal to $\alpha$ when $k'\neq k$ and to $dim_k(M)<\infty$ when $k'=k$ as there is no need for localization in this case. \\

A differential module $M$ is $0$-pure if and only if it is torsion-free. We have $t_0(M)=t(M)$ and $M/t(M)$ is torsion-free. In that case, using an {\it absolute localization}, one can find an {\it absolute parametrization} by embedding $M$ into a free module of rank $\alpha={\alpha}^n_q< m$, that is a parametrization by means of $\alpha$ arbitrary functions of $n$ variables. As a natural generalization, following Macaulay ([9]) while using a {\it relative localization} as in the previous corollaries, one can obtain ([18], [19]): \\

\noindent
{\bf COROLLARY 3.33}: When $M$ is $r$-pure, one can find a {\it relative parametrization} by means of $\alpha ={\alpha}^{n-r}_q$ arbitrary functions of $n-r$ variables, that is a parametrization by means of $\alpha $ functions which are {\it constrained} by a system of partial differential equations with no equation of class $n-r$ and full classes $n-r+1, ... ,n$ defining an $r$-pure parametrizing module $L$ with projective dimension equal to $r$, according to the corresponding Janet board of Section 2, both with an embedding $M\subseteq L$. (Compare to [1], p 494) (See [19], Section 4, for examples and counterexamples). \\    
      
\noindent
{\bf EXAMPLE 3.34}: With $n=3,m=1,k=\mathbb{Q}$, the module defined by the system $y_{33}=0, y_{23}=0, y_{22}=0, y_{13}-y_2=0$ is primary and 
$2$-pure. However, the module defined by the system $y_{33}=0, y_{23}=0, y_{22}=0, y_{13}=0$ is also of codimension $2$ but is not $2$-pure because we obtain by localizing/tensoring by $k({\chi}_1)$ the relations $y_{13}=0\Rightarrow {\chi}_1y_3=0\Rightarrow y_3=0$. Accordingly, the module defined by the subsystem $y_3=0, y_{22}=0$ is $2$-pure. We have nevertheless $\alpha=3-1=2$ and $par=\{y,y_2\}$ in both cases. Finally, the first system admits the general solution $y=f(x)=a(x^1)+{\partial}_1c(x^1)x^2+c(x^1)x^3$ which may not be easily related to the CK data (exercise). We let the reader work out the Spencer forms of the two systems and check Corollary 3.30 directly. Indeed, for the first system with $\alpha = 2$, setting $z^1=y, z^2=y_1, z^3=y_2, z^4=y_3 $, one obtains the relative parametrization:   \\
\[   (y,z) \leftrightarrow (z^1=y, z^2=y_1, z^3= z_1, z^4=z)  \]
coming from the localized equations ${\chi}_1z^1-z^2=0, {\chi}_1z^4-z^3=0$ with the following involutive differential constraint: \\
\[  \left\{  \begin{array}{rcl}
z_3=0 & , \hspace{5mm}&y_3-z=0 \\
z_2=0 & , \hspace{5mm} &y_2-z_1=0
\end{array}
\right.  \fbox {$ \begin{array}{rcl}
1 & 2 & 3 \\
1 & 2 & \bullet
\end{array} $ }  \]
coming from the other equations of class $2$ and $3$ by substitution. The elimination of $z$ provides of course the initial system of Macaulay and one obtains $z=g(x)=c(x^1)$ , that is a differential isomorphism $(y,z) \leftrightarrow (a(x^1), c(x^1))$ explaining this new concept.  \\

\noindent
{\bf EXAMPLE 3.35}: With $n=4, m=1, q=5, k=\mathbb{Q}$, we shall use modern methods in order to revisit the quite tricky example presented by Macaulay in ([9],$\S 85$,p 94), changing slightly the notations. The fifth order system defined by two equations of order $5$ and four equations of order $4$ as follows:  \\
\[    y_{23344}-y_{11333}=0 \hspace{1cm}, \hspace{1cm} y_{22444}-y_{13344}=0   \]
\[     y_{4444}=0, y_{3444}=0, y_{3334}=0, y_{3333}=0   \]
is neither formally integrable nor involutive. However, using one prolongation, we may obtain the following equivalent involutive system in which we have separated the equations of order $4$ and $5$ both with the various classes at order $5$, namely $1$ of class $4$, $5$ of class $3$, $6$ of class $2$ and $4$ of class $1$:  \\
\[  \left\{  \begin{array}{l}
y_{44444}=0   \\
y_{34444}=0,  y_{33444}=0, y_{33344}=0, y_{33334}=0, y_{33333}=0  \\
y_{24444}=0, y_{23444}=0, y_{23344}-y_{11333}=0, y_{23334}=0, y_{23333}=0, y_{22444}-y_{13344}=0  \\
y_{14444}=0, y_{13444}=0, y_{13334}=0, y_{13333}=0  \\
y_{4444}=0, y_{3444}=0, y_{3334}=0, y_{3333}=0
\end{array}   
\right.    \]
Such a system defines a differential module $M$ with $cd(M)=2$ and $\alpha=15-6=9$ because the classes $3$ and $4$ are full. This module is $2$-pure (exercise) but the interest of such an example is to examine the localized system with respect to $k({\chi}_1,{\chi}_2)=k({\chi}')=k'$ through the following involutive system of finite type with $5$ equations of order $4$ and {\it a new equation of order} $3$:  \\
\[  \left\{   \begin{array}{l}
y_{4444}=0 \\
y_{3444}=0 , \hspace{2mm}y_{3344}-\frac{({\chi}_1)^2}{{\chi}_2}y_{333}=0 ,\hspace{2mm} y_{3334}=0,\hspace{2mm} y_{3333}=0  \\
y_{444}-\frac{({\chi}_1)^3}{({\chi}_2)^3}y_{333}=0
\end{array}
\right.  \fbox{$\begin{array}{cc}
3  & 4  \\
3 & \bullet  \\
\bullet & \bullet
\end{array}$}  \]
Dividing by ${\chi}_2$ in a coherent way with the proof of the last corollary, we get for example:  \\
\[  \Phi\equiv y_{13444}=0, \Psi\equiv y_{23444}=0\Rightarrow d_2\Phi\equiv y_{123444}\equiv d_1\Psi \Rightarrow \Phi=({\chi}_1 / {\chi}_2)\Psi \]
We have thus $par_3=\{y,y_3,y_4, y_{33}, y_{34}, y_{44}, y_{333}, y_{334}, y_{344}\}$ providing $\alpha=9$ linearly independent modular equations/sections with a unique generating one with coefficients in $k[{\chi}_1,{\chi}_2]$, namely:  \\
\[    E\equiv ({\chi}_2)^3a^{333}+({\chi}_1)^3a^{444}+({\chi}_1{\chi}_2)^2a^{3344} = 0 \]
where {\it only one single determinant is needed} (we have chosen $({\chi}_2)^3$) in order to get polynomial coefficients ([18], Proposition 5.3.1). Indeed, applying the (Spencer) operators $d_3$ and $d_4$ in a convenient succession, one may obtain the $8$ additional modular equations/sections:  \\
\[  a^{344}=0, a^{334}=0, a^{44}=0, a^{34}=0, a^{33}=0, a^4=0, a^3=0, a^0=0  \]
a result not evident at first sight (exercise). We finally notice that $y_{333}$ is killed by $(d_3,d_4)\in max(k'[d_3,d_4])$, in a coherent way with ([18]). \\

We now explain and illustrate the way to use systems intead of modules for studying primary decompositions in the commutative framework as already explained. For simplicity, let us consider a primary decomposition with two components giving rise to a monomorphism $0\rightarrow M\rightarrow Q' \oplus Q''$ where $Q',Q''$ are primary modules, both with two epimorphisms $M\rightarrow Q' \rightarrow 0, M\rightarrow Q'' \rightarrow 0$, respectively induced by the localization morphisms $M\rightarrow M_{{\mathfrak{p}}{'}}, M\rightarrow M_{{\mathfrak{p}}{"}}$ when $M$ is pure (unmixed annihilator) with $ass(Q')=\{{\mathfrak{p}}'\}, ass(Q")=\{{\mathfrak{p}}"\}$ and $ass(M)=\{{\mathfrak{p}}{'},{\mathfrak{p}}{"}\}$. Setting $R'=hom_k(Q',k), R''=hom_k(Q'',k)$ and using the fact that $hom_k(D,k)$ is injective, we get an epimorphism $R' \oplus R'' \rightarrow R\rightarrow 0$ both with two monomorphisms $0\rightarrow R' \rightarrow R, 0\rightarrow R'' \rightarrow R$ proving that $R' , R'' , R' + R'' , R' \cap R'' $ are subsystems of $R$. The following proposition ([18], Prop 4.7), {\it not evident at first sight}, explains the aim of Macaulay ([9], end of \S 79, p 89) and allows to {\it use various subsystems for studying} $R$ {\it instead of decomposing} $M$.\\

\noindent
{\bf PROPOSITION 3.36}:  One has $ R = R' + R'' $. Accordingly $R \simeq R' \oplus R"$ or, equivalently, $ R' \cap R" =0$ if and only if 
${\mathfrak{p}}'+{\mathfrak{p}}"=A$.\\

\noindent
{\bf EXAMPLE 3.37}: With $n=3, m=1, q=2, k=\mathbb{Q}$, the module defined by the homogeneous involutive system $y_{33}=0, y_{23}-y_{13}=0, y_{22}-y_{12}=0$ is 2-pure. Setting $k({\chi}')=k({\chi}_1)$, the corresponding localized system  $y_{33}=0, y_{23}-{\chi}_1y_3=0, y_{22}-{\chi}_1y_2=0$ is again involutive with $par=\{y, y_2, y_3\}$. We obtain therefore a basis of sections $\{ f_1,f_2,f_3\}$ made by the three sections corresponding to $(1,0,0), (0,1,0), (0,0,1)$. We notice that ${\mathfrak{m}}"=(d_3,d_2-{\chi}_1)=ann(y_3)$ while ${\mathfrak{m}}'=(d_3,d_2)=ann(y_2-{\chi}_1y)$, each maximum ideal in $k({\chi}_1)[d_2,d_3]$ leading to a unique isotypical component of the socle of $M$. Denoting simply by $M$ the localized module and by $R$ the corresponding system with subsystems $R'$ provided by ${\mathfrak{m}}'$ and $R"$ provided by ${\mathfrak{m}}"$. It is then easy to obtain the primary decomposition (exercise) $\mathfrak{a}=({{\chi}_3}^2, {\chi}_3({\chi}_2-{\chi}_1),{\chi}_2({\chi}_2-{\chi}_1))=(({{\chi}_3}^2,{\chi}_2-{\chi}_1)\cap ({\chi}_3,{\chi}_2)={\mathfrak{q}}' \cap {\mathfrak{q}}" $ with $rad(\mathfrak{a})={\mathfrak{p}}' \cap {\mathfrak{p}}" $. After localization, ${\mathfrak{p}}' $ is replaced by ${\mathfrak{m}}' $ while ${\mathfrak{p}}"$ is replaced by ${\mathfrak{m}}"$ with $d_2-(d_2-{\chi}_1)={\chi}_1$ and thus ${\mathfrak{m}}' +{\mathfrak{m}}" =k({\chi}_1)[d_2,d_3]$. As for the systems, $R'$ is defined by $y_{33}=0, y_2-{\chi}_1 y=0$ while $R"$ is defined by $y_3=0,y_2=0$ and we have $R' \cap R"=0, R' \oplus R" \simeq R' + R" = R \Rightarrow dim(R)=dim(R')+dim(R")=2+1=3$. Finally, $R$ can be generated by the unique section $f=f_2+f_3$ because $d_2f + {\chi}_1 f= - f_1, d_3f= - f_1-{\chi}_1f_2$. The reader may study the involutive system $y_{33}-y_3=0, y_{23}-y_2=0, y_{22}-y_{12}=0, y_{13}-y_2=0$ providing a similar situation (exercise). \\

We are now ready for using the results of the second section on the Cartan-K\"{a}hler theorem. For such a purpose, we may write the solved equations in the symbolic form $y_{pri} - c^{par}_{pri} y_{par}=0$ with $c\in k$ and an implicit (finite) summation in order to obtain for the sections $f_{pri} - c^{par}_{pri}f_{par}=0$. Using the language of Macaulay, it follows that the so-called {\it modular equations} are $E\equiv f_{pri}a^{pri} + f_{par}a^{par}=0$ with eventually an infinite number of terms in the implicit summations. Substituting, we get at once $E\equiv f_{par}(a^{par} + c^{par}_{pri}a^{pri})=0$. Ordering the $y_{par}$ as we already did and using a basis $\{ (1,0,...), (0,1,0,...), (0,0,1,0,...), ... \}$ for the $f_{par}$, we may select the {\it parametric modular equations} $E^{par}\equiv a^{par} + c^{par}_{pri}a^{pri}=0$ and the same procedure could be used for the (finite type) localized system with $k'$ in place of $k$ and a finite number of such equations (Compare to (A)+(B) in [9], $\S 79$ or to (1)+(2) in [12]). \\

When a polynomial $P=a^{\mu}{\chi}^{\mu}\in k[\chi]$ of degree $q$ is multiplied by a monomial ${\chi}_{\nu}$ with $\mid \nu \mid=r$, we get ${\chi}_{\nu}P=a^{\mu}{\chi}_{\mu+\nu}$. Hence, if $0\leq \mid \mu\mid \leq q$, the "{\it shifted} " polynomial thus obtained is such that $r\leq \mid \mu+\nu\mid \leq q+r$ and the difference between the maximum degree and the minimum degree of the monomials involved is always equal to $q$ and thus fixed. This comment will allow to provide an example of a $k[\chi]$-linearly independent element of $k[\chi]^*=hom_k(k[\chi],k)$, namely a map $\varphi={\varphi}_n:k[\chi]\rightarrow k$ with ${\varphi}_0=id_k$ such that ${\chi}_{\mu}\rightarrow 1$ if ${\mu}_1=...={\mu}_n=p(p+3)/2$ for $p=0,1,2,...$ and $0$ otherwise. The reason for such a choice is that $p(p+3)/2-(p-1)(p+2)/2=p+1$ is strictly increasing with $p$. Indeed, if ${\mu}_1\leq d_1, ... ,{\mu}_n\leq d_n$ and $a_{(p_1,...,p_n)}\neq 0$, we may choose $p=max\{p_1,...,p_n\}$, set $d=p(p+3)/2$ and shift $P$ by $({\chi}_1)^{d-p_1}...({\chi}_n)^{d-p_n}$ in such a way that the contraction of $\varphi$ by this shifted operator is equal to $a_{(p_1,...,p_n)}\neq 0$, all the other lower terms beeing in what we shall call a "{\it zero zone}" of $\varphi$ of length $p+1>p$. Replacing ${\chi}_i$ by $d_i$ and {\it degree} by {\it order}, we may use the results of section 2 in order to split the CK-data into $m$ formal power series of $0$ (constants), $1,...,n$ variables that we shall call series of {\it type} $i$ for $i=0,1,...,n$. When a linear differential polynomial $P=P_ky^k=a^{\mu}_ky^k_{\mu}\in Dy\simeq D^m$ is given we may at once reduce it by using the previous reduction formulas in order to keep only the $y_{par}$ as we have only a finite number of them. We can thus decompose the reduction of $P$ into $m$ disjoint components, each one belonging to a series of a certain type $i=i(k)$ for $i=0,1,...,n$. For this, we shall set $y_{par}=\{y^k_{par}\mid 1\leq k \leq m\}$ with $y^k_{par}=\{y^k_{\mu}\mid {\mu}_{i+1}=...={\mu}_n=0,\forall {\beta}^{i+1}_1\leq k\leq {\beta}^i_1\}$ for $1\leq i\leq n-1$, $y^k_{par}=\{y^k_{\mu}\mid \forall {\beta}^n_1\leq k \leq m\}$ and $y^k_{par}=y^k\mid \forall 0\leq k \leq {\beta}^1_1\}$ as a way to define $i(k)=0$ for $1\leq k \leq {\beta}^1_1$, $i(k)=i$ for ${\beta}^{i+1}_1 \leq k \leq {\beta}^i_1$ and $i(k)=n$ for ${\beta}^n_1\leq k \leq m$. It just remains to compose each section $f^k_{par}=\{f^k_{\mu}\}$ with ${\varphi}_{i(k)}$ in order to have values equal to $0$ or $1$. Such a procedure can be extended {\it mutatis mutandis} by linearity to the parametric modular equations $E^{par}\equiv a^{par}+...=0$ when $D=k[d]$ in order to obtain $m$ formal power series with zero zones which cannot thus be killed by any operator. \\
However, as the following elementary example will show, this procedure cannot be applied to the variable coefficient case, namely when $K$ is used in place of $k$. Indeed, with $n=1, d_x=d, K=\mathbb{Q}(x)$ and $P=d^2-\frac{x}{3}$, if we contract $d^3P=d^5-\frac{x}{3}d^3-d^2$ with the series $f=(1,0,1,0,0,1,0,...)$ already defined when $n=1$, we get $1-1=0$ though $P$ does not kill $f$ because $d^2f=(1,0,0,1,0,0,...)\Rightarrow Pf=(1-\frac{x}{3},...)$ and the contraction of $P$ with $f$ is $1-\frac{x}{3}\neq 0$. WE SHALL ESCAPE FROM THIS DIFFICULTY BY MEANS OF A TRICK BASED ON A SYSTEMATIC USE OF THE SPENCER OPERATOR, shifting the series to the {\it left} ({\it decreasing ordering}), up to sign, instead of shifting the operator to the {\it right} ({\it increasing ordering}). For this, we notice that we want that the contraction of $P=a^{\mu}d_{\mu}$ where $\mid \mu \mid \leq q=ord(P)$ with $f$ should be zero, that is $a^{\mu}f_{\mu}=0\Rightarrow ({\partial}_ia^{\mu})f_{\mu}+a^{\mu}({\partial}_if_{\mu})=0,\forall i=1,...,n$. But $d_iP=a^{\mu}d_{\mu+1_i}+({\partial}_ia^{\mu})d_{\mu}$ must also contract to zero wih $f$ that is $a^{\mu}f_{\mu+1_i}+({\partial}_ia^{\mu})f_{\mu}=0$. Substracting, we obtain therefore the condition $a^{\mu}({\partial}_if_{\mu}-f_{\mu+1_i})=0$, that is $P$ must also contract to zero with the shift $d_if$ or even $d_{\nu}f$ of $f$ when $f$ is made with $0$ and $1$ only. Applying this computation to the above example, we get $-df=(0,1,0,0,1,0,...)\Rightarrow d^2f=(1,0,0,1,0,...)\Rightarrow -d^3f=(0,0,1,0,...)$ and the contraction with $P$ provides the leading coefficient $1\neq 0$ of $P$ like the contraction of $d^3P$ with $d^4f$, that is the same series can be used but in a quite different framework. We have therefore obtained the main result of this paper, in a coherent way with the finite dimensional case existing when the symbol $g_q$ of the defining system $R_q$ is finite type, that is when $g_{q+r}=0$ for a certain integer $r\geq 0$. Indeed, applying the $\delta$-sequence inductively to $g_{q+n+i}$ for $i=r-1,..., 0$ as in ([14], Proposition 6,p 87), it is known that $g_q$ is finite type {\it and} involutive if and only if $g_q=0$, that is to say $dim(R_q)=dim(R_{q-1})$:ÊÊ\\

\vspace{1cm}
\noindent
{\bf THEOREM 3.38}: If $M$ is a differential module over $D=K[d]$ defined by a first order involutive system in the $m$ unknowns $y^1,...,y^m$ with no zero order equation, the differential module $R=hom_K(M,K)$ may be generated over $D$ by a {\it finite basis of sections} containing $m$ generators. \\

In the general situation, counting the number of CK data, we have ${\alpha}^1_q+...+{\alpha}^n_q=dim(g_q)$ and $dim(R_q)=dim(g_q)+dim(R_{q-1})$. We obtain therefore the following result which is coherent with the number of unknowns in the Spencer form $R_{q+1}\subset J_1(R_q)$ (Compare to Theorem 2.3.1 in [8]):  \\

\noindent
{\bf COROLLARY 3.39}: If $M$ is a differential module over $D=K[d]$ defined by an involutive system $R_q\subset J_q(E)$, the differential module $R=hom_K(M,K)$ may be generated over $D$ by a {\it finite basis of sections} containing $dim(R_q)$ generators.  \\

\noindent
{\bf 4  CONCLUSION} \\

In 1916, Macaulay discovered new localization techniques for studying polynomial ideals after transforming them into systems of partial differential equations in one unknown. As a byproduct, he discovered the concept of formal integrability that will be studied later on successively by Riquier, Janet and Gr\"{o}bner for many unknowns but in a computational way. Using the intrinsic methods of the formal theory of systems of partial differential equations developped after 1960 by Spencer and coworkers, both with its application to the study of differential modules, we have been able to revisit the work of Macaulay and extend it from the constant to the variable coefficient case. In particular, the duality existing between differential modules and differential systems is crucially used in order to provide for the first time a link between the so-called inverse systems of Macaulay and the Cartan-K\"{a}hler theorem known for involutive systems, even for systems with coefficients in a given differential field. We hope that the many and sometimes tricky illustrating examples presented through this paper will become test examples for a future use of computer algebra.  \\

\noindent
{\bf REFERENCES}  \\

\noindent
[1] J.E. BJORK: Analytic D-Modules and Applications, Kluwer, 1993.\\
\noindent
[2] N. BOURBAKI: El\'{e}ments de Math\'{e}matiques, Alg\`{e}bre, Ch. 10, Alg\`{e}bre Homologique, Masson, Paris, 1980.   \\
\noindent
[3] E. CARTAN: Les Syst\`{e}mes Diff\'{e}rentiels Ext\'{e}rieurs et Leurs Applications G\'{e}om\'{e}triques, Hermann, Paris,1945.  \\
\noindent
[4] W. GR\"{O}BNER: \"{U}ber die Algebraischen Eigenschaften der Integrale von Linearen Differentialgleichungen mit Konstanten Koeffizienten, 
Monatsh. der Math., 47, 1939, 247-284.\\
\noindent
[5] M. JANET: Sur les Syst\`{e}mes aux D\'{e}riv\'{e}es Partielles, Journal de Math., 8(3),1920, 65-151.  \\
\noindent
[6] E. KAHLER: Einf\"{u}hrung in die Theorie der Systeme von Differentialgleichungen, Teubner, 1934. \\
\noindent
[7] M. KASHIWARA: Algebraic Study of Systems of Partial Differential Equations, M\'emoires de la Soci\'et\'e Math\'ematique de France 63, 1995, 
(Transl. from Japanese of his 1970 Master's Thesis).\\
\noindent
[8] E. KUNZ: Introduction to Commutative Algebra and Algebraic Geometry, Birkh\"{a}user, Boston, 1985.  \\
\noindent
[9] F. S. MACAULAY, The Algebraic Theory of Modular Systems, Cambridge Tracts 19, Cambridge University Press, London, 1916; Reprinted by Stechert-Hafner Service Agency, New York, 1964.\\
\noindent
[10] U. OBERST: Multidimensional Constant Linear Systems, Acta Appl. Math., 20, 1990, 1-175.   \\
\noindent
[11] U. OBERST: The Computation of Purity Filtrations over Commutative Noetherian Rings of Operators and their Applications to Behaviours, Multidim. Syst. Sign. Process. (MSSP), Springer, 2013.\\
http://dx.doi.org/10.1007/s11045-013-0253-4  \\
\noindent
[12] F. PIRAS: Some Remarks on the Inverse Systems of Polynomial Modules, Journal of Pure and Applied Algebra, 129,1998, 87-99. \\
\noindent
[13] J.-F. POMMARET: Systems of Partial Differential Equations and Lie Pseudogroups, Gordon and Breach, New York, 1978 (Russian translation by MIR, Moscow, 1983) \\
\noindent
[14] J.-F. POMMARET: Partial Differential Equations and Group Theory,New Perspectives for Applications, Mathematics and its Applications 293, Kluwer, 1994.\\
http://dx.doi.org/10.1007/978-94-017-2539-2   \\
\noindent
[15] J.-F. POMMARET: Partial Differential Control Theory, Kluwer, 2001, 957 pp.\\
\noindent
[16] J.-F. POMMARET: Algebraic Analysis of Control Systems Defined by Partial Differential Equations, in Advanced Topics in Control Systems Theory, Lecture Notes in Control and Information Sciences 311, Chapter 5, Springer, 2005, 155-223.\\
\noindent
[17] J.-F. POMMARET: Gr\"{o}bner Bases in Algebraic Analysis: New perspectives for applications, Radon Series Comp. Appl. Math 2, 1-21, de Gruyter, 2007.\\
\noindent
[18] J.-F. POMMARET: Macaulay Inverse Systems Revisited, Journal of Symbolic Computation, 46, 2011, 1049-1069.  \\
http://dx.doi.org/10.1016/j.jsc.2011.05.007    \\
\noindent
[19] J.-F. POMMARET: Relative Parametrization of Linear Multidimensional Systems, Multidim Syst Sign Process (MSSP), Springer, 2013. \\
http://dx.doi.org/10.1007/s11045-013-0265-0  \\
\noindent
[20] J.-F. POMMARET: The Mathematical Foundations of General Relativity Revisited, Journal of Modern Physics, 2013, 4,223-239.  \\
http://dx.doi.org/10.4236/jmp.2013.48A022  \\
\noindent
[21] J.-F. POMMARET: The Mathematical Foundations of Gauge Theory Revisited, Journal of Modern Physics, 2014, 5, 157-170.  \\
http://dx.doi.org/10.4236/jmp.2014.55026    \\
\noindent
[22] A. QUADRAT, D. ROBERTZ: A Constructive Study of the Module Structure of Rings of Partial Differential Operators, Acta Applicandae Mathematicae, 2014, (to appear). \\
http://hal-supelec.archives-ouvertes.fr/hal-00925533  \\
\noindent
[23] C. RIQUIER: Les Syst\`{e}mes d'Equations aux D\'{e}riv\'{e}es Partielles, Gauthiers-Villars, Paris, 1910.  \\
\noindent
[24] D. ROBERTZ: Formal Computational Methods for Control Theory, Dissertation, Mathematics Department, Aachen University, 20.06.2006.  \\
\noindent
[25] J.J. ROTMAN: An Introduction to Homological Algebra, Pure and Applied Mathematics, Academic Press, 1979.\\
\noindent
[26] J.-P. SCHNEIDERS: An Introduction to D-Modules, Bull. Soc. Roy. Sci. Li\`{e}ge, 63, 1994, 223-295.  \\
\noindent
[27] D.C. SPENCER: Overdetermined Systems of Partial Differential Equations, Bull. Amer. Math. Soc., 75, 1965, 1-114.\\

\end{document}